\documentclass[english,sort&compress]{elsarticle}
\usepackage[T1]{fontenc}
\usepackage[latin9]{inputenc}
\usepackage{geometry}
\usepackage{units}
\usepackage{amsmath}
\usepackage{amssymb}
\usepackage{cancel}
\usepackage{graphicx}

\makeatletter

\providecommand{\tabularnewline}{\\}

\usepackage{babel}

\makeatother

\usepackage{babel}
\begin{document}
\begin{frontmatter}
\title{Exact local conservation of energy in fully implicit PIC algorithms}
\author[lanl]{L. Chacón\corref{cor1}}
\ead{chacon@lanl.gov}
\author[lanl]{G. Chen}
\cortext[cor1]{Corresponding author}
\address[lanl]{Los Alamos National Laboratory, Los Alamos, NM 87545}
\begin{abstract}
We consider the issue of strict, fully discrete \emph{local} energy
conservation for a whole class of fully implicit local-charge- and
global-energy-conserving particle-in-cell (PIC) algorithms. Earlier
studies demonstrated these algorithms feature strict global energy
conservation. However, whether a local energy conservation theorem
exists (in which the local energy update is governed by a flux balance
equation at every mesh cell) for these schemes is unclear. In this
study, we show that a local energy conservation theorem indeed exists.
We begin our analysis with the 1D electrostatic PIC model without
orbit-averaging, and then generalize our conclusions to account for
orbit averaging, multiple dimensions, and electromagnetic models (Darwin).
In all cases, a temporally, spatially, and particle-discrete local
energy conservation theorem is shown to exist, proving that these
formulations (as originally proposed in the literature), in addition
to being locally charge conserving, are strictly locally energy conserving
as well. In contrast to earlier proofs of local conservation in the
literature \citep{xiao2017local}, which only considered continuum
time, our result is valid for the fully implicit time-discrete version
of all models, including important features such as orbit averaging.
We demonstrate the local-energy-conservation property numerically
with a paradigmatic numerical example. 
\end{abstract}
\end{frontmatter}

\section{Introduction}

Recently \citep{xiao2017local}, a proof of local energy conservation
in the time-continuum, spatially discrete context has been produced
for geometric particle-in-cell (PIC) methods (i.e., derived from a
suitable Lagrangian formulation). Since the proof in the reference
is for continuous time, it is unclear whether it survives temporal
discretization. Nevertheless, the approach outlines a template with
which to analyze other existing PIC algorithms, which while not derived
geometrically, may still feature strong conservation properties.

Over the last few years, there have been a series of studies \citep{chen2011energy,chen20141ddarwin,chacon2013curvilinear,chen2015multi,chacon2016curvilinear,chen2023implicit}
proposing locally charge-conserving, globally energy-conserving PIC
algorithms based on fully implicit time-centered discretizations.
These algorithms are remarkably flexible in that they allow the introduction
of orbit-averaging (critical to enforce charge conservation with particle
cell crossings),and curvilinear meshes without loss of conservation
properties. Moreover, they have been shown to be quite robust against
so-called aliasing (or finite-grid) instabilities~\citep{barnes2021finite}.

The study in Ref. \citep{xiao2017local} prompts the question as to
whether these fully implicit PIC algorithms feature a local energy
conservation theorem or not. In this study, we demonstrate that in
fact they all do so, at least in Cartesian geometry, and that the
property carries through to the temporally discrete setting, which
to our knowledge is a new contribution. This result strengthens the
already solid theoretical foundation of these fully implicit PIC algorithms,
and underpins their excellent performance in practice.

This paper is organized as follows. We begin in Sec. \ref{sec:darwin-cont}
with the derivation of the continuum local-energy conservation theorem
in the Vlasov-Darwin model (a low-frequency approximation to Maxwell's
equations \citep{krause2007unified}). We then begin in Sec. \ref{sec:1d-es}
the exploration of the discrete subtleties with the simplest model,
a one-dimensional electrostatic PIC. From there, we build up in complexity
first with orbit averaging (Sec. \ref{sec:orbit-averaging}), then
in dimensionality (Sec. \ref{sec:multi-d}), and later by including
electromagnetic physics via the Darwin model (Sec. \ref{sec:em-md}).
Our proofs leverage a key theorem proved in Ref. \citep{xiao2017local}
that demonstrates the flux form of a remaining numerical error term,
key to cast the conservation equation in a flux form. In Sec. \ref{sec:numerics},
we demonstrate numerically with a one-dimensional electrostatic example
that the local energy conservation theorem is indeed enforced by the
scheme, and that the numerical error term can be quantified and is
in fact quite small. Finally, we conclude in Sec. \ref{sec:conclusions}.

\section{Continuum local energy conservation equation in the Vlasov-Darwin
model}

\label{sec:darwin-cont}

We begin this study by deriving the continuum local energy conservation
equation for the Vlasov-Darwin model in multiple dimensions \citep{krause2007unified,chen2015multi},
which is the ultimate target of this study. This model is comprised
of the Vlasov equation for multiple species $\alpha$: 
\[
\partial_{t}f_{\alpha}+\mathbf{v}\cdot\nabla f_{\alpha}+\frac{q_{\alpha}}{m_{\alpha}}(\mathbf{E}+\mathbf{v}\times\mathbf{B})\cdot\nabla_{v}f_{\alpha}=0,
\]
coupled with the low-frequency Ampere's equation: 
\begin{equation}
-\partial_{t}\nabla\phi+\mathbf{j}=\nabla\times\nabla\times\mathbf{A}\label{eq:ampere-cont}
\end{equation}
and its divergence: 
\begin{equation}
\partial_{t}\nabla^{2}\phi=\nabla\cdot\mathbf{j}.\label{eq:charge-cons}
\end{equation}
Here, 
\[
\mathbf{j}=\sum_{\alpha}q_{\alpha}\int d\mathbf{v}\mathbf{v}f_{\alpha}
\]
is the current density. The vector potential $\mathbf{A}$ gives the
magnetic field as: 
\begin{equation}
\mathbf{B}=\nabla\times\mathbf{A},\label{eq:B-field-def}
\end{equation}
and is in the Coulomb gauge, and therefore solenoidal: 
\[
\nabla\cdot\mathbf{A}=0.
\]
The electric field is given by: 
\begin{equation}
\mathbf{E}=-\nabla\phi-\partial_{t}\mathbf{A}.\label{eq:E-field-def}
\end{equation}
Owing to the Coulomb gauge, the longitudinal component is electrostatic,
while the vector potential time derivative is purely inductive.

Taking the second moment of Vlasov's equation, $m\int dv\,v^{2}/2[...]$,
we find the local energy conservation equation: 
\begin{equation}
\partial_{t}e_{k}+\nabla\cdot\boldsymbol{\Gamma}=\mathbf{E}\cdot\mathbf{j},\label{eq:Ek-local-cons}
\end{equation}
with the kinetic energy density $e_{k}$ and flux $\boldsymbol{\Gamma}$
given by: 
\[
e_{k}=m\int dv\frac{v^{2}}{2}f\,\,,\,\,\boldsymbol{\Gamma}=m\int dv\frac{v^{2}}{2}\mathbf{v}f.
\]
Dotting Ampere's equation with $\nabla\phi$, we find: 
\[
-\partial_{t}\frac{(\nabla\phi)^{2}}{2}+\mathbf{j}\cdot\nabla\phi=(\nabla\times\nabla\times\mathbf{A})\cdot\nabla\phi,
\]
where $\mathcal{E}=(\nabla\phi)^{2}/2$ is the electrostatic energy.
Using the definition of the electric field, Eq. \ref{eq:E-field-def},
and Ampere's equation, we find: 
\[
-\partial_{t}\mathcal{E}-\mathbf{j}\cdot\mathbf{E}=(\nabla\times\nabla\times\mathbf{A})\cdot\nabla\phi+(\nabla\times\nabla\times\mathbf{A}+\partial_{t}\nabla\phi)\cdot\partial_{t}\mathbf{A},
\]
which can be rewritten using the Coulomb gauge as: 
\[
\partial_{t}\mathcal{E}+\mathbf{j}\cdot\mathbf{E}=-\nabla\cdot(\phi\nabla\times\nabla\times\mathbf{A}+\partial_{t}\phi\cdot\partial_{t}\mathbf{A})-(\nabla\times\nabla\times\mathbf{A})\cdot\partial_{t}\mathbf{A}.
\]
The measure of the local magnetic energy in terms of the vector potential
can be written as: 
\[
\mathcal{M}=\frac{1}{2}(\nabla\times\mathbf{A})^{2}=\frac{1}{2}\nabla\mathbf{A}:(\nabla\mathbf{A}-\nabla\mathbf{A}^{\dagger})=\frac{1}{4}(\nabla\mathbf{A}-\nabla\mathbf{A}^{\dagger}):(\nabla\mathbf{A}-\nabla\mathbf{A}^{\dagger}).
\]
Also, we note that: 
\begin{equation}
\nabla\times\nabla\times\mathbf{A}=-\nabla\cdot(\nabla\mathbf{A}-\nabla\mathbf{A}^{\dagger})=-\nabla^{2}\mathbf{A},\label{eq:veclap_def}
\end{equation}
where the dagger symbol stands for transpose. The term $\nabla\cdot\nabla\mathbf{A}^{\dagger}=\nabla(\nabla\cdot\mathbf{A)}=0$
due to the Coulomb gauge, but it is advantageous to keep it for clarity
in the derivation. Then, we can write: 
\[
-(\nabla\times\nabla\times\mathbf{A})\cdot\partial_{t}\mathbf{A}=\nabla\cdot\left((\nabla\mathbf{A}-\nabla\mathbf{A}^{\dagger})\cdot\partial_{t}\mathbf{A}\right)-(\nabla\mathbf{A}-\nabla\mathbf{A}^{\dagger}):\partial_{t}\nabla\mathbf{A}=\nabla\cdot\left((\nabla\mathbf{A}-\nabla\mathbf{A}^{\dagger})\cdot\partial_{t}\mathbf{A}\right)-\partial_{t}\mathcal{M}.
\]
There results the EM field local conservation theorem: 
\begin{equation}
\partial_{t}\mathcal{M}+\partial_{t}\mathcal{E}+\mathbf{j}\cdot\mathbf{E}=\nabla\cdot\boldsymbol{\Gamma}_{EM},\label{eq:EM-local-cons}
\end{equation}
with the field-energy flux given by: 
\begin{equation}
\boldsymbol{\Gamma}_{EM}=\phi\nabla^{2}\mathbf{A}-\partial_{t}\phi\cdot\partial_{t}\mathbf{A}+(\nabla\mathbf{A}-\nabla\mathbf{A}^{\dagger})\cdot\partial_{t}\mathbf{A}.\label{eq:em-flux}
\end{equation}
The field-energy flux can be reinterpreted as a measure of the Poynting
flux by noting that: 
\begin{eqnarray*}
\phi\nabla^{2}\mathbf{A} & = & -\phi\nabla\times\mathbf{B}=-\nabla\times(\phi\mathbf{B})+\nabla\phi\times\mathbf{B},\\
(\nabla\mathbf{A}-\nabla\mathbf{A}^{\dagger})\cdot\partial_{t}\mathbf{A} & = & \partial_{t}\mathbf{A}\times\nabla\times\mathbf{A}=\partial_{t}\mathbf{A}\times\mathbf{B},
\end{eqnarray*}
to find: 
\[
\nabla\cdot\boldsymbol{\Gamma}_{EM}=\nabla\cdot(-\mathbf{E}\times\mathbf{B}-\partial_{t}\phi\partial_{t}\mathbf{A}).
\]
The $\nabla\cdot(\partial_{t}\phi\partial_{t}\mathbf{A})$ term is
a correction to the standard Poynting's flux stemming from the Darwin
approximation.

Combining Eq. \ref{eq:EM-local-cons} with the kinetic-energy local
conservation law, Eq. \ref{eq:Ek-local-cons}, we finally find the
total energy conservation law sought: 
\begin{equation}
\partial_{t}\left(e_{k}+\mathcal{E}+\mathcal{M}\right)+\nabla\cdot(\boldsymbol{\Gamma}+\mathbf{E}\times\mathbf{B}+\partial_{t}\phi\partial_{t}\mathbf{A})=0.\label{eq:local_e_cons_continuum}
\end{equation}
Global energy conservation follows by integrating over the whole domain
with appropriate boundary conditions.

\section{Local energy conservation theorem in a fully discrete 1D electrostatic
PIC model}

\label{sec:1d-es}

We begin the discrete analysis with the simplest fully implicit 1D
electrostatic PIC model. We will later generalize the results to multiple
dimensions and electromagnetic models. The time-discrete formulation
(without orbit averaging) of a fully implicit, Ampere-based 1D electrostatic
model reads \citep{chen2011energy}: 
\begin{eqnarray}
\frac{x_{p}^{n+1}-x_{p}^{n}}{\Delta t} & = & v_{p}^{n+1/2},\label{eq:pos}\\
\frac{v_{p}^{n+1}-v_{p}^{n}}{\Delta t} & = & \frac{q_{p}}{m_{p}}E_{p}^{n+1/2},\label{eq:eom}\\
E_{p}^{n+1/2} & = & \sum_{i}E_{i+1/2}^{n+1/2}S_{m-1}(x_{i+1/2}-x_{p}^{n+1/2})\label{eq:E_scatter}\\
\frac{E_{i+1/2}^{n+1}-E_{i+1/2}^{n}}{\Delta t}+j_{i+1/2}^{n+1/2} & = & \langle j\rangle,\label{eq:ampere}\\
j_{i+1/2}^{n+1/2} & = & \frac{1}{\Delta x}\sum_{p}q_{p}v_{p}^{n+1/2}S_{m-1}(x_{i+1/2}-x_{p}^{n+1/2}),\label{eq:j_gather}
\end{eqnarray}
where $n+1/2$ quantities are simple averages of integer time levels,
and $\langle j\rangle$ is the average current in a periodic system
(else, zero). Here, $m$ is the choice of B-spline interpolation order,
which is typically of second order.

We begin by defining the cell-centered local kinetic energy density
moment as: 
\begin{equation}
e_{k,i}=\frac{1}{\Delta x}\sum_{p}m_{p}\frac{v_{p}^{2}}{2}S_{l}(x_{i}-x_{p}).\label{eq:local-kin-en-1}
\end{equation}
As has been standard in earlier studies \citep{chen2011energy,chen20141ddarwin,chacon2013curvilinear,chen2015multi,chacon2016curvilinear,chen2023implicit},
we consider spline orders $l=1,2$. Local energy conservation laws
can be found for either choice. Equation \ref{eq:local-kin-en-1}
satisfies that the kinetic energy density spatial integral equates
the total particle kinetic energy in the system: 
\[
\sum_{i}\Delta xe_{k,i}=\sum_{p}m_{p}\frac{v_{p}^{2}}{2}
\]
owing to the interpolation kernels being a partition of unity. We
also consider a complementary local field energy measure: 
\[
\mathcal{E}_{i}=\frac{E_{i+1/2}^{2}+E_{i-1/2}^{2}}{4},
\]
which we will use to define a local field-energy conservation law.

\subsection{Kinetic energy local conservation law}

We begin with the discrete temporal derivative of the energy density
in Eq. \ref{eq:local-kin-en-1} as: 
\begin{equation}
\frac{e_{k,i}^{n+1}-e_{k,i}^{n}}{\Delta t}=\frac{1}{\Delta x\Delta t}\sum_{p}m_{p}\left[\frac{(v_{p}^{n+1})^{2}}{2}S_{l}(x_{i}-x_{p}^{n+1})-\frac{(v_{p}^{n})^{2}}{2}S_{l}(x_{i}-x_{p}^{n})\right].\label{eq:dEk_dt}
\end{equation}
For a B-spline up to second order and \emph{within a given cell} $i$,
we can exactly write by Taylor expansion: 
\begin{eqnarray}
S_{l}(x_{i}-x_{p}^{n+1}) & = & S_{l}(x_{i}-x_{p}^{n+1/2})-\left.\frac{\partial S_{l}}{\partial x_{p}}\right|_{i}^{n+1/2}(x_{p}^{n+1}-x_{p}^{n+1/2})+\left.\frac{\partial^{2}S_{l}}{\partial x_{p}^{2}}\right|_{i}^{n+1/2}\frac{(x_{p}^{n+1}-x_{p}^{n+1/2})^{2}}{2},\label{eq:taylor-np}\\
S_{l}(x_{i}-x_{p}^{n}) & = & S_{l}(x_{i}-x_{p}^{n+1/2})-\left.\frac{\partial S_{l}}{\partial x_{p}}\right|_{i}^{n+1/2}(x_{p}^{n}-x_{p}^{n+1/2})+\left.\frac{\partial^{2}S_{l}}{\partial x_{p}^{2}}\right|_{i}^{n+1/2}\frac{(x_{p}^{n}-x_{p}^{n+1/2})^{2}}{2}.\label{eq:taylor-n}
\end{eqnarray}
Introducing these into Eq. \ref{eq:dEk_dt} and manipulating (\ref{app:1d-es-derivation})
gives the following kinetic-energy \emph{local} conservation law:
\begin{equation}
\frac{e_{k,i}^{n+1}-e_{k,i}^{n}}{\Delta t}+\frac{\Gamma_{i+1/2}^{n+1/2}-\Gamma_{i-1/2}^{n+1/2}}{\Delta x}=\frac{1}{\Delta x}\sum_{p}q_{p}E_{p}^{n+1/2}v_{p}^{n+1/2}\tilde{S}_{l}(x_{i}-x_{p}^{n+1/2}),\label{eq:cons-e_k-final}
\end{equation}
where the energy fluxes are defined as: 
\[
\Gamma_{i+1/2}^{n+1/2}=\frac{1}{\Delta x}\sum_{p}S_{l-1}(x_{i+1/2}-x_{p}^{n+1/2})v_{p}^{n+1/2}e_{k,p}^{n+1/2},
\]
with 
\[
e_{k,p}^{n+1/2}=\frac{m_{p}}{2}\left[\frac{(v_{p}^{n+1})^{2}}{2}+\frac{(v_{p}^{n})^{2}}{2}\right]
\]
a time-centered particle kinetic energy, and: 
\[
\tilde{S}_{l}(x_{i}-x_{p}^{n+1/2})=S_{l}(x_{i}-x_{p}^{n+1/2})+\left.\frac{\partial^{2}S_{l}}{\partial x_{p}^{2}}\right|_{i}^{n+1/2}\frac{\Delta t^{2}}{8}\left(v_{p}^{n+1/2}\right)^{2},
\]
which remains a partition of unity, since: 
\begin{equation}
\sum_{i}\tilde{S}_{l}(x_{i}-x_{p}^{n+1/2})=1,\label{eq:partition-of-unity-2}
\end{equation}
because $\sum_{i}\left.\frac{\partial^{2}S_{l}}{\partial x_{p}^{2}}\right|_{i}^{n+1/2}=0$.

The source in Eq. \ref{eq:cons-e_k-final} approximates the $E.j$
term at the cell center $i$. Note that Eq. \ref{eq:cons-e_k-final}
still leads to the expected global conservation statement, since in
a periodic domain: 
\begin{eqnarray*}
\sum_{i}\Delta x\frac{\Gamma_{i+1/2}^{n+1/2}-\Gamma_{i-1/2}^{n+1/2}}{\Delta x} & = & 0,
\end{eqnarray*}
and, using Eq. \ref{eq:partition-of-unity-2}: 
\begin{equation}
\sum_{i}\Delta x\frac{1}{\Delta x}\sum_{p}q_{p}E_{p}^{n+1/2}v_{p}^{n+1/2}\tilde{S}_{l}(x_{i}-x_{p}^{n+1/2})=\sum_{p}q_{p}E_{p}^{n+1/2}v_{p}^{n+1/2}=\sum_{i}\Delta xj_{i+1/2}^{n+1/2}E_{i+1/2}^{n+1/2},\label{eq:E.j-result}
\end{equation}
and therefore: 
\[
\sum_{i}\Delta x\frac{e_{k,i}^{n+1}-e_{k,i}^{n}}{\Delta t}=\sum_{i}\Delta xj_{i+1/2}^{n+1/2}E_{i+1/2}^{n+1/2},
\]
from which the total energy conservation theorem follows \citep{chen2011energy}.
In Eq. \ref{eq:E.j-result}, the first equality follows from the partition
of unity of the B-splines and the fact that the particle quantities
are independent of the mesh point $i$, and the second equality follows
from the definition of the field scatter to the mesh (Eq.~\ref{eq:E_scatter}).

\subsection{Electric field energy temporal update equation}

We define the E-field energy at the $i$-cell as the average of the
field energy at faces: 
\[
\mathcal{E}_{i}^{n+1}=\frac{(E_{i+1/2}^{n+1})^{2}+(E_{i-1/2}^{n+1})^{2}}{4}.
\]
It follows that the local temporal rate of change of the E-field energy
is given by: 
\begin{eqnarray}
\frac{\mathcal{E}_{i}^{n+1}-\mathcal{E}_{i}^{n}}{\Delta t} & = & \frac{(E_{i+1/2}^{n+1})^{2}+(E_{i-1/2}^{n+1})^{2}-(E_{i+1/2}^{n})^{2}-(E_{i-1/2}^{n})^{2}}{4\Delta t}\nonumber \\
 & = & \frac{(E_{i+1/2}^{n+1}+E_{i+1/2}^{n})}{4}\frac{(E_{i+1/2}^{n+1}-E_{i+1/2}^{n})}{\Delta t}+\frac{(E_{i-1/2}^{n+1}+E_{i-1/2}^{n})}{4}\frac{(E_{i-1/2}^{n+1}-E_{i-1/2}^{n})}{\Delta t}\nonumber \\
 & = & -\frac{E_{i+1/2}^{n+1/2}(j_{i+1/2}^{n+1/2}-\langle j\rangle)+E_{i-1/2}^{n+1/2}(j_{i-1/2}^{n+1/2}-\langle j\rangle)}{2}\nonumber \\
 & = & -\frac{(E_{i+1/2}^{n+1/2}j_{i+1/2}^{n+1/2}+E_{i-1/2}^{n+1/2}j_{i-1/2}^{n+1/2})}{2}+\frac{(E_{i+1/2}^{n+1/2}+E_{i-1/2}^{n+1/2})}{2}\langle j\rangle.\label{eq:e-field-cons-law}
\end{eqnarray}
The latter term is a source term due to the possible presence of an
average current, and is a source of energy into the system. This term
can be written in conservative form by using the electrostatic potential
$\phi$, since: 
\begin{equation}
\frac{(E_{i+1/2}^{n+1/2}+E_{i-1/2}^{n+1/2})}{2}\langle j\rangle=-\frac{\left(\phi_{i+1}^{n+1/2}+\phi_{i}^{n+1/2}\right)-\left(\phi_{i-1}^{n+1/2}+\phi_{i}^{n+1/2}\right)}{2\Delta x}\langle j\rangle=\frac{\Gamma_{E,i+1/2}^{n+1/2}-\Gamma_{E,i-1/2}^{n+1/2}}{\Delta x},
\end{equation}
where: 
\begin{equation}
\Gamma_{E,i+1/2}^{n+1/2}=-\langle j\rangle\frac{\phi_{i+1}^{n+1/2}+\phi_{i}^{n+1/2}}{2}.\label{eq:avg_current_cont}
\end{equation}

\subsection{Total-energy local conservation law}

Adding Eqs. \ref{eq:cons-e_k-final} and \ref{eq:e-field-cons-law},
we find: 
\begin{eqnarray}
\frac{e_{k,i}^{n+1}-e_{k,i}^{n}}{\Delta t}+\frac{\mathcal{E}_{i}^{n+1}-\mathcal{E}_{i}^{n}}{\Delta t} & = & -\frac{\Gamma_{i+1/2}^{n+1/2}-\Gamma_{E,i+1/2}^{n+1/2}-\Gamma_{i-1/2}^{n+1/2}+\Gamma_{E,i-1/2}^{n+1/2}}{\Delta x}\nonumber \\
 & + & \frac{1}{\Delta x}\sum_{p}q_{p}v_{p}^{n+1/2}\tilde{S}_{l}(x_{i}-x_{p}^{n+1/2})\sum_{j}E_{j+1/2}^{n+1/2}S_{m-1}(x_{j+1/2}-x_{p}^{n+1/2})\nonumber \\
 & - & \frac{1}{\Delta x}\sum_{p}q_{p}v_{p}^{n+1/2}\frac{\left[E_{i+1/2}^{n+1/2}S_{m-1}(x_{i+1/2}-x_{p}^{n+1/2})+E_{i-1/2}^{n+1/2}S_{m-1}(x_{i-1/2}-x_{p}^{n+1/2})\right]}{2}.\label{eq:cons-e_k-1}
\end{eqnarray}
The source terms approximate the same continuum term ($\mathbf{E}\cdot\mathbf{j})$,
and therefore asymptotically offset one another with increasing mesh/particle
resolution (as will be shown numerically later in this study). More
importantly, the difference can be written in conservation form. Following
\citep{xiao2017local}, and defining: 
\begin{align*}
g_{j+1/2,p}^{n+1/2} & =v_{p}^{n+1/2}E_{j+1/2}^{n+1/2}S_{m-1}(x_{j+1/2}-x_{p}^{n+1/2}),\\
F_{p}^{n+1/2} & =\sum_{j}g_{j+1/2,p}^{n+1/2}=v_{p}^{n+1/2}E_{p}^{n+1/2},
\end{align*}
we can rewrite the source term in Eq. \ref{eq:cons-e_k-1} as: 
\begin{eqnarray*}
 &  & \frac{1}{\Delta x}\sum_{p}q_{p}v_{p}^{n+1/2}\tilde{S}_{l}(x_{i}-x_{p}^{n+1/2})\sum_{j}E_{j+1/2}^{n+1/2}S_{m-1}(x_{j+1/2}-x_{p}^{n+1/2})\\
 & - & \frac{1}{\Delta x}\sum_{p}q_{p}v_{p}^{n+1/2}\frac{\left[E_{i+1/2}^{n+1/2}S_{m-1}(x_{i+1/2}-x_{p}^{n+1/2})+E_{i-1/2}^{n+1/2}S_{m-1}(x_{i-1/2}-x_{p}^{n+1/2})\right]}{2}\\
 & = & \frac{1}{\Delta x}\sum_{p}q_{p}\left[F_{p}^{n+1/2}\tilde{S}_{l}(x_{i}-x_{p}^{n+1/2})-\frac{g_{i+1/2,p}^{n+1/2}+g_{i-1/2,p}^{n+1/2}}{2}\right].
\end{eqnarray*}
Note that this trivially satisfies the global energy conservation
theorem, since by definition: 
\begin{eqnarray*}
\sum_{i}\left\{ \sum_{p}q_{p}\left[F_{p}^{n+1/2}\tilde{S}_{l}(x_{i}-x_{p}^{n+1/2})-\frac{g_{i+1/2,p}^{n+1/2}+g_{i-1/2,p}^{n+1/2}}{2}\right]\right\}  & = & 0.
\end{eqnarray*}
Moreover, the scalar expression in curly brackets is a discrete, sum-free
(per the above equation), and locally non-zero (it is only finite
in a stencil around cell $i$) scalar field, and, as proven in Ref.
\citep{xiao2017local}, it can be written in a discrete divergence
form, i.e.: 
\[
\frac{1}{\Delta x}\sum_{p}q_{p}\left[F_{p}^{n+1/2}\tilde{S}_{l}(x_{i}-x_{p}^{n+1/2})-\frac{g_{i+1/2,p}^{n+1/2}+g_{i-1/2,p}^{n+1/2}}{2}\right]=\frac{\Gamma_{K,i+1/2}-\Gamma_{K,i-1/2}}{\Delta x}.
\]
and is therefore locally conservative. This last equation defines
the numerical error flux $\Gamma_{K,i+1/2}$ from its source. This
numerical flux will be quantified later in this study.

The final total-energy local conservation theorem finally reads: 
\[
\frac{e_{k,i}^{n+1}-e_{k,i}^{n}}{\Delta t}+\frac{\mathcal{E}_{i}^{n+1}-\mathcal{E}_{i}^{n}}{\Delta t}+\frac{\Gamma_{i+1/2}^{n+1/2}-\Gamma_{E,i+1/2}^{n+1/2}-\Gamma_{i-1/2}^{n+1/2}+\Gamma_{E,i-1/2}^{n+1/2}}{\Delta x}=\frac{\Gamma_{K,i+1/2}^{n+1/2}-\Gamma_{K,i-1/2}^{n+1/2}}{\Delta x}.
\]
This result is the specialization of Eq. \ref{eq:local_e_cons_continuum}
for the 1D electrostatic model. However, it assumed no cell crossings,
which is not the case in general. To account for particle crossings,
we need to introduce orbit averaging. We consider this extension next.

\section{Generalization of local-energy conservation equation to include orbit-averaging}

\label{sec:orbit-averaging}

The time-discrete formulation of the fully implicit 1D electrostatic
PIC model with orbit averaging~\citep{chen2011energy} (assuming
for simplicity $\langle j\rangle=0$ in this section) reads: 
\begin{eqnarray}
\frac{x_{p}^{\nu+1}-x_{p}^{\nu}}{\Delta\tau^{\nu}} & = & v_{p}^{\nu+1/2},\label{eq:pos-1}\\
\frac{v_{p}^{\nu+1}-v_{p}^{\nu}}{\Delta\tau^{\nu}} & = & \frac{q_{p}}{m_{p}}E_{p}^{\nu+1/2},\label{eq:eom-1}\\
E_{p}^{\nu+1/2} & = & \sum_{i}E_{i+1/2}^{n+1/2}S_{m-1}(x_{i+1/2}-x_{p}^{\nu+1/2})\label{eq:E_scatter-1}\\
\frac{E_{i+1/2}^{n+1}-E_{i+1/2}^{n}}{\Delta t}+\bar{j}_{i+1/2}^{n+1/2} & = & 0,\label{eq:ampere-1}\\
\bar{j}_{i+1/2}^{n+1/2} & = & \frac{1}{\Delta x\Delta t}\sum_{p}q_{p}\sum_{\nu}\Delta\tau^{\nu}v_{p}^{\nu+1/2}S_{m-1}(x_{i+1/2}-x_{p}^{\nu+1/2}),\label{eq:j_gather-1}
\end{eqnarray}
We consider the corresponding orbit-averaged kinetic-energy update
from $n$ to $n+1$ (from Eq. \ref{eq:local-kin-en-1}): 
\begin{equation}
\frac{e_{k,i}^{n+1}-e_{k,i}^{n}}{\Delta t}=\frac{1}{\Delta x\Delta t}\sum_{p}\sum_{\nu\in(n,n+1)}m_{p}\left[\frac{(v_{p}^{\nu+1})^{2}}{2}S_{l}(x_{i}-x_{p}^{\nu+1})-\frac{(v_{p}^{\nu})^{2}}{2}S_{l}(x_{i}-x_{p}^{\nu})\right].\label{eq:dEk_dt-1}
\end{equation}
Following the same procedure as in the previous section, we find:
\begin{eqnarray}
\frac{e_{k,i}^{n+1}-e_{k,i}^{n}}{\Delta t} & = & \frac{1}{\Delta x\Delta t}\sum_{p}\sum_{\nu\in(n,n+1)}\Delta\tau^{\nu}q_{p}E_{p}^{\nu+1/2}v_{p}^{\nu+1/2}\tilde{S}_{l}(x_{i}-x_{p}^{\nu+1/2})\nonumber \\
 & - & \frac{1}{\Delta x\Delta t}\sum_{p}\sum_{\nu\in(n,n+1)}\Delta\tau^{\nu}\frac{S_{l-1}(x_{i+1/2}-x_{p}^{\nu+1/2})-S_{l-1}(x_{i-1/2}-x_{p}^{\nu+1/2})}{\Delta x}v_{p}^{\nu+1/2}e_{k,p}^{\nu+1/2}\nonumber \\
 & = & S^{n+1/2}-\frac{\Gamma_{i+1/2}^{n+1/2}-\Gamma_{i-1/2}^{n+1/2}}{\Delta x},\label{eq:local-kin-energy-cons-orb_avg}
\end{eqnarray}
where, as before: 
\[
e_{k,p}^{\nu+1/2}=\frac{m_{p}}{2}\left[\frac{(v_{p}^{\nu+1})^{2}}{2}+\frac{(v_{p}^{\nu})^{2}}{2}\right].
\]
This equation is identical to Eq. \ref{eq:cons-e_k-final}, but with
orbit-averaged source and fluxes: 
\begin{eqnarray*}
S^{n+1/2} & = & \frac{1}{\Delta x\Delta t}\sum_{p}\sum_{\nu\in(n,n+1)}\Delta\tau^{\nu}q_{p}E_{p}^{\nu+1/2}v_{p}^{\nu+1/2}\tilde{S}_{l}(x_{i}-x_{p}^{\nu+1/2}),\\
\Gamma_{i+1/2}^{n+1/2} & = & \frac{1}{\Delta x\Delta t}\sum_{p}\sum_{\nu\in(n,n+1)}\Delta\tau^{\nu}S_{l-1}(x_{i+1/2}-x_{p}^{\nu+1/2})v_{p}^{\nu+1/2}e_{k,p}^{\nu+1/2}.
\end{eqnarray*}
As before, summing over the whole mesh gives the expected result:
\[
\sum_{i}\Delta x\frac{e_{k,i}^{n+1}-e_{k,i}^{n}}{\Delta t}=\frac{1}{\Delta t}\sum_{p}\sum_{\nu\in(n,n+1)}\Delta\tau^{\nu}q_{p}E_{p}^{\nu+1/2}v_{p}^{\nu+1/2}=\sum_{i}\Delta xE_{i+1/2}^{n+1/2}\bar{j}_{i+1/2}^{n+1/2}.
\]

Regarding the local rate of change of the electric field energy, as
before we have: 
\begin{eqnarray}
\frac{\mathcal{E}_{i}^{n+1}-\mathcal{E}_{i}^{n}}{\Delta t} & = & \frac{(E_{i+1/2}^{n+1})^{2}+(E_{i-1/2}^{n+1})^{2}-(E_{i+1/2}^{n})^{2}-(E_{i-1/2}^{n})^{2}}{4\Delta t}\nonumber \\
 & = & \frac{(E_{i+1/2}^{n+1}+E_{i+1/2}^{n})}{4}\frac{(E_{i+1/2}^{n+1}-E_{i+1/2}^{n})}{\Delta t}+\frac{(E_{i-1/2}^{n+1}+E_{i-1/2}^{n})}{4}\frac{(E_{i-1/2}^{n+1}-E_{i-1/2}^{n})}{\Delta t}\nonumber \\
 & = & -\frac{(E_{i+1/2}^{n+1/2}\bar{j}_{i+1/2}^{n+1/2}+E_{i-1/2}^{n+1/2}\bar{j}_{i-1/2}^{n+1/2})}{2}\label{eq:e-field-energy-1d}
\end{eqnarray}
where now the current density is orbit averaged per Eq. \ref{eq:ampere-1}.
Combining this result with the kinetic energy conservation law (Eq.
\ref{eq:local-kin-energy-cons-orb_avg}), we find: 
\begin{eqnarray}
\frac{e_{k,i}^{n+1}-e_{k,i}^{n}}{\Delta t}+\frac{\mathcal{E}_{i}^{n+1}-\mathcal{E}_{i}^{n}}{\Delta t} & = & -\frac{\Gamma_{i+1/2}^{n+1/2}-\Gamma_{i-1/2}^{n+1/2}}{\Delta x}\nonumber \\
 & + & \frac{1}{\Delta x\Delta t}\sum_{p}q_{p}\sum_{\nu\in(n,n+1)}\Delta\tau^{\nu}E_{p}^{\nu+1/2}v_{p}^{\nu+1/2}\tilde{S}_{l}(x_{i}-x_{p}^{\nu+1/2})\nonumber \\
 & - & \frac{1}{2\Delta x\Delta t}\sum_{p}q_{p}\sum_{\nu\in(n,n+1)}\Delta\tau^{\nu}v_{p}^{\nu+1/2}E_{i+1/2}^{n+1/2}S_{m-1}(x_{i+1/2}-x_{p}^{\nu+1/2})\nonumber \\
 & - & \frac{1}{2\Delta x\Delta t}\sum_{p}q_{p}\sum_{\nu\in(n,n+1)}\Delta\tau^{\nu}v_{p}^{\nu+1/2}E_{i-1/2}^{n+1/2}S_{m-1}(x_{i-1/2}-x_{p}^{\nu+1/2}).\label{eq:cons-e_k-orbit_avg}
\end{eqnarray}
Following the same recipe, defining: 
\begin{align*}
g_{j+1/2,p}^{\nu+1/2} & =E_{j+1/2}^{n+1/2}v_{p}^{\nu+1/2}S_{m-1}(x_{j+1/2}-x_{p}^{\nu+1/2}),\\
F_{p}^{\nu+1/2} & =\sum_{j}g_{j+1/2,p}^{\nu+1/2}=E_{p}^{\nu+1/2}v_{p}^{\nu+1/2},
\end{align*}
the source can then be written as: 
\[
\frac{1}{\Delta x\Delta t}\sum_{p}\sum_{\nu\in(n,n+1)}\Delta\tau^{\nu}q_{p}\left[F_{p}^{\nu+1/2}\tilde{S}_{l}(x_{i}-x_{p}^{\nu+1/2})-\frac{g_{i+1/2,p}^{\nu+1/2}+g_{i-1/2,p}^{\nu+1/2}}{2}\right].
\]
As before, this quantity is locally non-zero, discrete, and sums to
zero: 
\[
\sum_{i}\left\{ \frac{1}{\Delta x\Delta t}\sum_{p}\sum_{\nu\in(n,n+1)}\Delta\tau^{\nu}q_{p}\left[F_{p}^{\nu+1/2}\tilde{S}_{l}(x_{i}-x_{p}^{\nu+1/2})-\frac{g_{i+1/2,p}^{\nu+1/2}+g_{i-1/2,p}^{\nu+1/2}}{2}\right]\right\} =0,
\]
and therefore can be written in in divergence form: 
\[
\frac{1}{\Delta x}\left\{ \frac{1}{\Delta t}\sum_{p}\sum_{\nu\in(n,n+1)}\Delta\tau^{\nu}q_{p}\left[F_{p}^{\nu+1/2}\tilde{S}_{l}(x_{i}-x_{p}^{\nu+1/2})-\frac{g_{i+1/2,p}^{\nu+1/2}+g_{i-1/2,p}^{\nu+1/2}}{2}\right]\right\} =\frac{\Gamma_{K,i+1/2}^{\nu+1/2}-\Gamma_{K,i-1/2}^{\nu+1/2}}{\Delta x},
\]
where again $\Gamma_{K}$ is a numerical flux, to be quantified later
in this study. This result leads to the \emph{local} energy conservation
law sought: 
\[
\frac{e_{k,i}^{n+1}-e_{k,i}^{n}}{\Delta t}+\frac{\mathcal{E}_{i}^{n+1}-\mathcal{E}_{i}^{n}}{\Delta t}+\frac{\Gamma_{i+1/2}^{n+1/2}-\Gamma_{i-1/2}^{n+1/2}}{\Delta x}=\frac{\Gamma_{K,i+1/2}^{n+1/2}-\Gamma_{K,i-1/2}^{n+1/2}}{\Delta x}.
\]
This result is identical to the no-cell-crossing result in the previous
section, but with orbit-averaged fluxes. We consider the generalization
to multiple dimensions next.

\section{Generalization to multiple dimensions}

\label{sec:multi-d}

In multiple dimensions, the fully implicit, orbit-averaged electrostatic
PIC algorithm reads \citep{chen2015multi}: 
\begin{eqnarray}
\frac{\mathbf{x}_{p}^{\nu+1}-\mathbf{x}_{p}^{\nu}}{\Delta\tau^{\nu}} & = & \mathbf{v}_{p}^{\nu+1/2},\label{eq:pos-md}\\
\frac{\mathbf{v}_{p}^{\nu+1}-\mathbf{v}_{p}^{\nu}}{\Delta\tau} & = & \frac{q_{p}}{m_{p}}\mathbf{E}_{p}^{\nu+1/2},\label{eq:eom-md}\\
\mathbf{E}_{p}^{\nu+1/2} & = & \sum_{g\in i,j,k}\mathbf{E}_{g+1/2}^{n+1/2}\cdot\bar{\bar{\mathbf{S}}}(\mathbf{x}_{g+1/2}-\mathbf{x}_{p}^{\nu+1/2}),\label{eq:E_scatter-md}\\
\mathbf{E}_{g+1/2}^{n+1/2} & = & -\nabla_{h}\phi_{g}^{n+1/2},\label{eq:e-field-def}\\
\frac{\nabla_{h}^{2}\phi_{g}^{n+1}-\nabla_{h}^{2}\phi_{g}^{n}}{\Delta t} & = & \nabla_{h}\cdot\mathbf{\bar{j}}_{g+1/2}^{n+1/2},\label{eq:ampere-md}\\
\mathbf{\bar{j}}_{g+1/2}^{n+1/2} & = & \frac{1}{\Delta_{h}\Delta t}\sum_{p}q_{p}\sum_{\nu\in(n,n+1)}\Delta\tau^{\nu}\mathbf{v}_{p}^{\nu+1/2}\cdot\bar{\bar{\mathbf{S}}}(\mathbf{x}_{g+1/2}-\mathbf{x}_{p}^{\nu+1/2}).\label{eq:j_gather-md}
\end{eqnarray}
In these equations, $g\in(i,j,k)$, $\phi_{g}=\phi_{i,j,k}$ is cell
centered, $\mathbf{E}_{g+1/2}=(E_{x,i+1/2,j,k},E_{y,i,j+1/2,k},E_{z,i,j,k+1/2})$
(and same with $\bar{\mathbf{j}}_{g+1/2}$) is face centered, and
we define: 
\begin{eqnarray}
\nabla_{h}\phi_{g} & = & \frac{\phi_{i+1,j,k}-\phi_{i,j,k}}{\Delta x}\mathbf{i}+\frac{\phi_{i,j+1,k}-\phi_{i,j,k}}{\Delta y}\mathbf{j}+\frac{\phi_{i,j,k+1}-\phi_{i,j,k}}{\Delta z}\mathbf{k},\\
\nabla_{h}\cdot\boldsymbol{\Gamma}_{g+1/2}^{n+1/2} & = & \frac{\Gamma_{x,i+1/2,j,k}-\Gamma_{x,i-1/2,j,k}}{\Delta x}+\frac{\Gamma_{y,i,j+1/2,k}-\Gamma_{y,i,j-1/2,k}}{\Delta y}+\frac{\Gamma_{z,i,j,k+1/2}-\Gamma_{z,i,j,k-1/2}}{\Delta z},
\end{eqnarray}
where $\mathbf{i},\mathbf{j},\mathbf{k}$ are unit vectors in the
$x$, $y$, $z$ directions, respectively. Also, $\Delta_{h}=\Delta x\Delta y\Delta z$,
and the face-centered shape-function dyad at $g+1/2$ is defined as:
\begin{eqnarray}
\bar{\bar{\mathbf{S}}}(\mathbf{x}_{g+1/2}-\mathbf{x}_{p}^{\nu+\nicefrac{1}{2}}) & = & \mathbf{i}\otimes\mathbf{i}\,\,S_{1}(x_{i+\nicefrac{1}{2}}-x_{p}^{\nu+\nicefrac{1}{2}})\mathcal{\mathbb{S}}_{22,jk}^{\nu+\nicefrac{1}{2}}(y_{p},z_{p})\nonumber \\
 & + & \mathbf{j}\otimes\mathbf{j}\,\,S_{1}(y_{j+\nicefrac{1}{2}}-y_{p}^{\nu+\nicefrac{1}{2}})\mathcal{\mathbb{S}}_{22,ik}^{\nu+\nicefrac{1}{2}}(z_{p},x_{p})\label{eq:S-doublebar}\\
 & + & \mathbf{k}\otimes\mathbf{k}\,\,S_{1}(z_{k+\nicefrac{1}{2}}-z_{p}^{\nu+\nicefrac{1}{2}})\mathcal{\mathbb{S}}_{22,ij}^{\nu+\nicefrac{1}{2}}(x_{p},y_{p}),\nonumber 
\end{eqnarray}
where $\otimes$ denotes tensor product, $S_{1}(x_{i+\nicefrac{1}{2}}-x_{p}^{\nu+\nicefrac{1}{2}})$
denotes linear B-spline shape function with $x_{i+\nicefrac{1}{2}}$
at cell faces, and 
\begin{align}
\mathcal{\mathbb{S}}_{22,jk}^{\nu+\nicefrac{1}{2}}(y,z)\equiv\frac{1}{3} & \left[S_{2}(y_{j}-y^{\nu+1})S_{2}(z_{k}-z^{\nu+1})+\frac{S_{2}(y_{j}-y^{\nu})S_{2}(z_{k}-z^{\nu+1})}{2}\right.\nonumber \\
 & \left.\:\;+\frac{S_{2}(y_{j}-y^{\nu+1})S_{2}(z_{k}-z^{\nu})}{2}+S_{2}(y_{j}-y^{\nu})S_{2}(z_{k}-z^{\nu})\right]\label{eq:S_22}
\end{align}
is defined with the explicit purpose of enforcing charge conservation
\citep{chen2015multi}. In the last equation, $S_{2}(y_{j}-y^{\nu+1})$
is the quadratic B-spline shape function, with $y_{j}$ denoting a
cell center (in $y$). The charge density and other cell-centered
scalars are gathered as: 
\[
\rho_{g}^{\nu}=\frac{1}{\Delta_{h}}\sum_{p}q_{p}S_{2}(\mathbf{x}_{g}-\mathbf{x}_{p}^{\nu}),
\]
where the vector argument in the shape function indicates tensor product,
e.g.: 
\[
S_{2}(\mathbf{x}_{g}-\mathbf{x}_{p}^{\nu})=S_{2}(x_{i}-x_{p}^{\nu})S_{2}(y_{j}-y_{p}^{\nu})S_{2}(z_{k}-z_{p}^{\nu}).
\]

\subsection{Kinetic-energy conservation law}

We consider the orbit-averaged kinetic-energy update from $n$ to
$n+1$ (from Eq. \ref{eq:local-kin-en-1}): 
\begin{equation}
\frac{e_{k,g}^{n+1}-e_{k,g}^{n}}{\Delta t}=\frac{1}{\Delta_{h}\Delta t}\sum_{p}\sum_{\nu\in(n,n+1)}m_{p}\left[\frac{(v_{p}^{\nu+1})^{2}}{2}S_{l}(\mathbf{x}_{g}-\mathbf{x}_{p}^{\nu+1})-\frac{(v_{p}^{\nu})^{2}}{2}S_{l}(\mathbf{x}_{g}-\mathbf{x}_{p}^{\nu})\right].\label{eq:dEk_dt-multi-D}
\end{equation}
Taylor-expanding the shape functions and manipulating (\ref{app:=00003D000020multi-D-es-derivation})
gives: 
\begin{eqnarray*}
\frac{e_{k,g}^{n+1}-e_{k,g}^{n}}{\Delta t}+\nabla_{h}\cdot\boldsymbol{\Gamma}_{g}^{n+1/2} & = & S_{g}^{n+1/2},
\end{eqnarray*}
where the source $S_{g}^{n+1/2}$ and the components of the flux $\boldsymbol{\Gamma}_{g+1/2}^{n+1/2}$
are given as: 
\begin{eqnarray}
S_{g}^{n+1/2} & = & \frac{1}{\Delta_{h}\Delta t}\sum_{p}\sum_{\nu\in(n,n+1)}\Delta\tau^{\nu}q_{p}\mathbf{E}_{p}^{\nu+1/2}\cdot\mathbf{v}_{p}^{\nu+1/2}\tilde{\tilde{S}}_{l}(\mathbf{x}_{g}-\mathbf{x}_{p}^{\nu+1/2}),\label{eq:md-source}\\
\Gamma_{x,i+1/2,j,k}^{\nu+1/2} & = & \frac{1}{\Delta_{h}\Delta t}\sum_{p}\sum_{\nu\in(n,n+1)}\Delta\tau^{\nu}S_{1}(x_{i+1/2}-x_{p}^{\nu+1/2})S_{2}(y_{j}-y_{p}^{\nu+1/2})S_{2}(z_{k}-z_{p}^{\nu+1/2})v_{x,p}^{\nu+1/2}e_{k,p}^{\nu+1/2},\label{eq:md-flux-x}\\
\Gamma_{y,i,j+1/2,k}^{\nu+1/2} & = & \frac{1}{\Delta_{h}\Delta t}\sum_{p}\sum_{\nu\in(n,n+1)}\Delta\tau^{\nu}S_{2}(x_{i}-x_{p}^{\nu+1/2})S_{1}(y_{j+1/2}-y_{p}^{\nu+1/2})S_{2}(z_{k}-z_{p}^{\nu+1/2})v_{y,p}^{\nu+1/2}e_{k,p}^{\nu+1/2},\label{eq:md-flux-y}\\
\Gamma_{z,i,j,k+1/2}^{\nu+1/2} & = & \frac{1}{\Delta_{h}\Delta t}\sum_{p}\sum_{\nu\in(n,n+1)}\Delta\tau^{\nu}S_{2}(x_{i}-x_{p}^{\nu+1/2})S_{2}(y_{j}-y_{p}^{\nu+1/2})S_{1}(z_{k+1/2}-z_{p}^{\nu+1/2})v_{z,p}^{\nu+1/2}e_{k,p}^{\nu+1/2}.\label{eq:md-flux-z}
\end{eqnarray}
Here: 
\[
\tilde{\tilde{S}}_{l}(\mathbf{x}_{g}-\mathbf{x}_{p}^{\nu+1/2})=S_{l}(\mathbf{x}_{g}-\mathbf{x}_{p}^{\nu+1/2})+\frac{(\Delta\tau^{\nu})^{2}}{8}\mathbf{v}_{p}^{\nu+1/2}\cdot(\left.\nabla_{x_{p}}\nabla_{x_{p}}S_{l}\right|_{g}^{\nu+1/2})\cdot\mathbf{v}_{p}^{\nu+1/2},
\]
and, as before: 
\[
e_{k,p}^{\nu+1/2}=\frac{m_{p}}{2}\left[\frac{(v_{p}^{\nu+1})^{2}}{2}+\frac{(v_{p}^{\nu})^{2}}{2}\right].
\]
Also as before, summing over the whole mesh gives the global conservation
theorem for kinetic energy: 
\[
\sum_{g}\Delta_{h}\frac{e_{k,g}^{n+1}-e_{k,g}^{n}}{\Delta t}=\frac{1}{\Delta t}\sum_{p}\sum_{\nu\in(n,n+1)}\Delta\tau^{\nu}q_{p}\mathbf{E}_{p}^{\nu+1/2}\cdot\mathbf{v}_{p}^{\nu+1/2}\sum_{g}\tilde{\tilde{S}}_{l}(\mathbf{x}_{g}-\mathbf{x}_{p}^{\nu+1/2})=\sum_{g\in(i,j,k)}\Delta_{h}\mathbf{E}_{g+1/2}^{n+1/2}\cdot\bar{\mathbf{j}}_{g+1/2}^{n+1/2},
\]
where we have used that: 
\begin{equation}
\sum_{g}\tilde{\tilde{S}}_{l}(\mathbf{x}_{g}-\mathbf{x}_{p}^{\nu+1/2})=\sum_{g}S_{l}(\mathbf{x}_{g}-\mathbf{x}_{p}^{\nu+1/2})+\frac{(\Delta\tau^{\nu})^{2}}{8}\mathbf{v}_{p}^{\nu+1/2}\cdot\nabla_{x_{p}}\nabla_{x_{p}}\sum_{g}\left.S_{l}\right|_{g}^{\nu+1/2}\cdot\mathbf{v}_{p}^{\nu+1/2}=1,\label{eq:partition_unity_S_tld_tld}
\end{equation}
since: 
\[
\sum_{g}S_{l}(\mathbf{x}_{g}-\mathbf{x}_{p}^{\nu+1/2})=1\,\,;\,\,\nabla_{x_{p}}\nabla_{x_{p}}\sum_{g}S_{l}(\mathbf{x}_{g}-\mathbf{x}_{p}^{\nu+1/2})=0.
\]

\subsection{Electrostatic energy local conservation law}

Defining the electrostatic energy as: 
\[
\mathcal{E}_{g}=\frac{(E_{x,i+1/2,j,k})^{2}+(E_{x,i-1/2,j,k})^{2}+(E_{y,i,j+1/2,k})^{2}+(E_{y,i,j-1/2,k})^{2}+(E_{z,i,j,k+1/2})^{2}+(E_{z,i,j,k-1/2})^{2}}{4},
\]
we find the multi-dimensional generalization of the electrostatic
energy conservation law in Eq. \ref{eq:e-field-energy-1d} as: 
\begin{eqnarray}
\left.\frac{\mathcal{E}_{g}^{n+1}-\mathcal{E}_{g}^{n}}{\Delta t}\right|_{g=i,j,k} & = & \frac{1}{2}\left[\mathbf{E}_{g+1/2}^{n+1/2}\cdot\frac{\mathbf{E}_{g+1/2}^{n+1}-\mathbf{E}_{g+1/2}^{n}}{\Delta t}+\mathbf{E}_{g-1/2}^{n+1/2}\cdot\frac{\mathbf{E}_{g-1/2}^{n+1}-\mathbf{E}_{g-1/2}^{n}}{\Delta t}\right]\nonumber \\
 & = & -\frac{1}{2}\left[\frac{\mathbf{E}_{g+1/2}^{n+1}-\mathbf{E}_{g+1/2}^{n}}{\Delta t}\cdot\nabla_{h}\phi_{g}^{n+1/2}+\frac{\mathbf{E}_{g-1/2}^{n+1}-\mathbf{E}_{g-1/2}^{n}}{\Delta t}\cdot\nabla_{h}\phi_{g-1}^{n+1/2}\right].\label{eq:e-field-energy-3d}
\end{eqnarray}
Now, we note the following \emph{exact} discrete identity (used often
from now on; proof by expansion): 
\begin{equation}
\frac{1}{2}\left[A_{i+1/2}\frac{\phi_{i+1}-\phi_{i}}{\Delta x}+A_{i-1/2}\frac{\phi_{i}-\phi_{i-1}}{\Delta x}\right]=\frac{1}{\Delta x}\left[\frac{\phi_{i+1}+\phi_{i}}{2}A_{i+1/2}-\frac{\phi_{i-1}+\phi_{i}}{2}A_{i-1/2}\right]-\phi_{i}\frac{A_{i+1/2}-A_{i-1/2}}{\Delta x},\label{eq:exact-chain-rule-div}
\end{equation}
from which it follows that (using Eq. \ref{eq:ampere-md}): 
\begin{eqnarray*}
 &  & \frac{1}{2}\left[\frac{\mathbf{E}_{g+1/2}^{n+1}-\mathbf{E}_{g+1/2}^{n}}{\Delta t}\cdot\nabla_{h}\phi_{g}^{n+1/2}+\frac{\mathbf{E}_{g-1/2}^{n+1}-\mathbf{E}_{g-1/2}^{n}}{\Delta t}\cdot\nabla_{h}\phi_{g-1}^{n+1/2}\right]\\
 & = & \nabla_{h}\cdot\left(\phi_{g+1/2}^{n+1/2}\frac{\mathbf{E}_{g+1/2}^{n+1}-\mathbf{E}_{g+1/2}^{n}}{\Delta t}\right)-\phi_{g}^{n+1/2}\nabla_{h}\cdot\left(\frac{\mathbf{E}_{g+1/2}^{n+1}-\mathbf{E}_{g+1/2}^{n}}{\Delta t}\right).\\
 & = & \nabla_{h}\cdot\left(\phi_{g+1/2}^{n+1/2}\frac{\mathbf{E}_{g+1/2}^{n+1}-\mathbf{E}_{g+1/2}^{n}}{\Delta t}\right)+\phi_{g}^{n+1/2}\nabla_{h}\cdot\bar{\mathbf{j}}_{g+1/2}^{n+1/2}.
\end{eqnarray*}
Here, $\phi_{g+1/2}$ is the average of cell-centered values to the
corresponding face. Similarly, from the same identity we can exactly
write: 
\[
\phi_{g}^{n+1/2}\nabla_{h}\cdot\bar{\mathbf{j}}_{g+1/2}^{n+1/2}=\nabla_{h}\cdot\left(\phi_{g+1/2}^{n+1/2}\bar{\mathbf{j}}_{g+1/2}^{n+1/2}\right)-\frac{1}{2}\left[\bar{\mathbf{j}}_{g+1/2}^{n+1/2}\cdot\nabla_{h}\phi_{g}^{n+1/2}+\bar{\mathbf{j}}_{g-1/2}^{n+1/2}\cdot\nabla_{h}\phi_{g-1}^{n+1/2}\right].
\]
Introducing these identities into Eq. \ref{eq:e-field-energy-3d},
there results: 
\[
\frac{\mathcal{E}_{g}^{n+1}-\mathcal{E}_{g}^{n}}{\Delta t}=\nabla_{h}\cdot\left(-\phi_{g+1/2}^{n+1/2}\left(\frac{\mathbf{E}_{g+1/2}^{n+1}-\mathbf{E}_{g+1/2}^{n}}{\Delta t}+\bar{\mathbf{j}}_{g+1/2}^{n+1/2}\right)\right)-\frac{1}{2}\left[\bar{\mathbf{j}}_{g+1/2}^{n+1/2}\cdot\mathbf{E}_{g+1/2}^{n+1/2}+\bar{\mathbf{j}}_{g-1/2}^{n+1/2}\cdot\mathbf{E}_{g-1/2}^{n+1}\right],
\]
which is the local electrostatic energy conservation equation sought.
We can finally write: 
\begin{eqnarray*}
\frac{\mathcal{E}_{g}^{n+1}-\mathcal{E}_{g}^{n}}{\Delta t} & = & \nabla_{h}\cdot\boldsymbol{\Gamma}_{E,g+1/2}^{n+1/2}-\frac{\mathbf{E}_{g+1/2}^{n+1/2}\cdot\bar{\mathbf{j}}_{g+1/2}^{n+1/2}+\mathbf{E}_{g-1/2}^{n+1/2}\cdot\bar{\mathbf{j}}_{g-1/2}^{n+1/2}}{2}\\
 & = & \nabla_{h}\cdot\boldsymbol{\Gamma}_{E,g+1/2}^{n+1/2}\\
 & - & \frac{1}{2\Delta_{h}\Delta t}\sum_{p}q_{p}\sum_{\nu\in(n,n+1)}\Delta\tau^{\nu}\mathbf{v}_{p}^{\nu+1/2}\cdot\bar{\bar{\mathbf{S}}}(\mathbf{x}_{g+1/2}-\mathbf{x}_{p}^{\nu+1/2})\cdot\mathbf{E}_{g+1/2}^{n+1/2}\\
 & - & \frac{1}{2\Delta_{h}\Delta t}\sum_{p}q_{p}\sum_{\nu\in(n,n+1)}\Delta\tau^{\nu}\mathbf{v}_{p}^{\nu+1/2}\cdot\bar{\bar{\mathbf{S}}}(\mathbf{x}_{g-1/2}-\mathbf{x}_{p}^{\nu+1/2})\cdot\mathbf{E}_{g-1/2}^{n+1/2},
\end{eqnarray*}
with: 
\begin{equation}
\boldsymbol{\Gamma}_{E,g+1/2}^{n+1/2}=-\phi_{g+1/2}^{n+1/2}\left(\frac{\mathbf{E}_{g+1/2}^{n+1}-\mathbf{E}_{g+1/2}^{n}}{\Delta t}+\bar{\mathbf{j}}_{g+1/2}^{n+1/2}\right),\label{eq:E-energy-flux}
\end{equation}
which is the generalization of Eq. \ref{eq:avg_current_cont}.

\subsection{Total-energy local conservation law}

Putting all together, we find: 
\begin{eqnarray*}
\frac{e_{k,g}^{n+1}-e_{k,g}^{n}}{\Delta t}+\frac{\mathcal{E}_{g}^{n+1}-\mathcal{E}_{g}^{n}}{\Delta t} & = & \nabla_{h}\cdot\left(-\boldsymbol{\Gamma}_{g+1/2}^{n+1/2}+\boldsymbol{\Gamma}_{E,g+1/2}^{n+1/2}\right)\\
 & + & \frac{1}{\Delta_{h}\Delta t}\sum_{p}q_{p}\sum_{\nu\in(n,n+1)}\Delta\tau^{\nu}\mathbf{E}_{p}^{\nu+1/2}\cdot\mathbf{v}_{p}^{\nu+1/2}\tilde{\tilde{S}}_{l}(\mathbf{x}_{g}-\mathbf{x}_{p}^{\nu+1/2})\\
 & - & \frac{1}{2\Delta_{h}\Delta t}\sum_{p}q_{p}\sum_{\nu\in(n,n+1)}\Delta\tau^{\nu}\mathbf{v}_{p}^{\nu+1/2}\cdot\bar{\bar{\mathbf{S}}}(\mathbf{x}_{g+1/2}-\mathbf{x}_{p}^{\nu+1/2})\cdot\mathbf{E}_{g+1/2}^{n+1/2}\\
 & - & \frac{1}{2\Delta_{h}\Delta t}\sum_{p}q_{p}\sum_{\nu\in(n,n+1)}\Delta\tau^{\nu}\mathbf{v}_{p}^{\nu+1/2}\cdot\bar{\bar{\mathbf{S}}}(\mathbf{x}_{g-1/2}-\mathbf{x}_{p}^{\nu+1/2})\cdot\mathbf{E}_{g-1/2}^{n+1/2}.
\end{eqnarray*}
As before, the global energy conservation theorem is immediately recovered,
since: 
\[
\sum_{g}\left[\mathbf{E}_{p}^{\nu+1/2}\cdot\mathbf{v}_{p}^{\nu+1/2}\tilde{\tilde{S}}_{l}(\mathbf{x}_{g}-\mathbf{x}_{p}^{\nu+1/2})-\mathbf{v}_{p}^{\nu+1/2}\cdot\frac{\bar{\bar{\mathbf{S}}}(\mathbf{x}_{g+1/2}-\mathbf{x}_{p}^{\nu+1/2})\cdot\mathbf{E}_{g+1/2}^{n+1/2}+\bar{\bar{\mathbf{S}}}(\mathbf{x}_{g-1/2}-\mathbf{x}_{p}^{\nu+1/2})\cdot\mathbf{E}_{g-1/2}^{n+1/2}}{2}\right]
\]
vanishes by Eqs. \ref{eq:E_scatter-md} and \ref{eq:partition_unity_S_tld_tld}.
Defining: 
\[
G_{g+1/2,p}^{\nu+1/2}=\mathbf{E}_{g+1/2}^{n+1/2}\cdot\bar{\bar{\mathbf{S}}}(\mathbf{x}_{g+1/2}-\mathbf{x}_{p}^{\nu+1/2})\cdot\mathbf{v}_{p}^{\nu+1/2}\,\,;\,\,F_{p}^{\nu+1/2}=\sum_{g}G_{g+1/2,p}^{\nu+1/2}=\mathbf{E}_{p}^{\nu+1/2}\cdot\mathbf{v}_{p}^{\nu+1/2},
\]
we can write: 
\begin{eqnarray*}
 &  & \mathbf{E}_{p}^{\nu+1/2}\cdot\mathbf{v}_{p}^{\nu+1/2}\tilde{\tilde{S}}_{l}(\mathbf{x}_{g}-\mathbf{x}_{p}^{\nu+1/2})-\mathbf{v}_{p}^{\nu+1/2}\cdot\frac{\bar{\bar{\mathbf{S}}}(\mathbf{x}_{g+1/2}-\mathbf{x}_{p}^{\nu+1/2})\cdot\mathbf{E}_{g+1/2}^{n+1/2}+\bar{\bar{\mathbf{S}}}(\mathbf{x}_{g-1/2}-\mathbf{x}_{p}^{\nu+1/2})\cdot\mathbf{E}_{g-1/2}^{n+1/2}}{2}\\
 & = & \left[F_{p}^{\nu+1/2}\tilde{\tilde{S}}_{l}(\mathbf{x}_{g}-\mathbf{x}_{p}^{\nu+1/2})-\frac{G_{g+1/2,p}^{\nu+1/2}+G_{g-1/2,p}^{\nu+1/2}}{2}\right].
\end{eqnarray*}
Again, it can be readily checked that the scalar in square brackets
is locally non-zero, discrete, and sums to zero: 
\[
\sum_{g}\left\{ \frac{1}{\Delta t}\sum_{p}q_{p}\sum_{\nu\in(n,n+1)}\Delta\tau^{\nu}\left[F_{p}^{\nu+1/2}\tilde{\tilde{S}}_{l}(\mathbf{x}_{g}-\mathbf{x}_{p}^{\nu+1/2})-\frac{G_{g+1/2,p}^{\nu+1/2}+G_{g-1/2,p}^{\nu+1/2}}{2}\right]\right\} =0,
\]
and therefore by the theorem in \citep{xiao2017local} the term in
curly brackets can be written in divergence form: 
\[
\frac{1}{\Delta_{h}}\frac{1}{\Delta t}\sum_{p}q_{p}\sum_{\nu\in(n,n+1)}\Delta\tau^{\nu}\left[F_{p}^{\nu+1/2}\tilde{\tilde{S}}_{l}(\mathbf{x}_{g}-\mathbf{x}_{p}^{\nu+1/2})-\frac{G_{g+1/2,p}^{\nu+1/2}+G_{g-1/2,p}^{\nu+1/2}}{2}\right]=\nabla_{h}\cdot\mathbf{\Gamma}_{K,g+1/2}^{\nu+1/2}.
\]
The total-energy local conservation law is now complete: 
\[
\frac{e_{k,g}^{n+1}-e_{k,g}^{n}}{\Delta t}+\frac{\mathcal{E}_{g}^{n+1}-\mathcal{E}_{g}^{n}}{\Delta t}+\nabla_{h}\cdot\left(\boldsymbol{\Gamma}_{g+1/2}^{n+1/2}-\boldsymbol{\Gamma}_{E,g+1/2}^{n+1/2}\right)=\nabla_{h}\cdot\boldsymbol{\Gamma}_{K,g+1/2}^{n+1/2}.
\]

\section{Generalization to EM (Darwin)}

\label{sec:em-md}

The discrete multidimensional Darwin-PIC model reads \citep{chen2015multi}:
\begin{eqnarray}
\frac{\mathbf{x}_{p}^{\nu+1}-\mathbf{x}_{p}^{\nu}}{\Delta\tau^{\nu}} & = & \mathbf{v}_{p}^{\nu+1/2},\label{eq:pos-md-1}\\
\frac{\mathbf{v}_{p}^{\nu+1}-\mathbf{v}_{p}^{\nu}}{\Delta\tau} & = & \frac{q_{p}}{m_{p}}\left(\mathbf{E}_{p}^{\nu+1/2}+\mathbf{v}_{p}^{\nu+1/2}\times\mathbf{B}_{p}^{\nu+1/2}\right),\label{eq:eom-md-1}\\
\mathbf{E}_{p}^{\nu+1/2} & = & \sum_{g\in i,j,k}\mathbf{E}_{g+1/2}^{n+1/2}\cdot\bar{\bar{\mathbf{S}}}(\mathbf{x}_{g+1/2}-\mathbf{x}_{p}^{\nu+1/2}),\label{eq:E_scatter-md-1}\\
\frac{\nabla_{h}^{2}\phi_{g}^{n+1}-\nabla_{h}^{2}\phi_{g}^{n}}{\Delta t} & = & \nabla_{h}\cdot\mathbf{\bar{j}}_{g+1/2}^{n+1/2},\label{eq:charge-md-em}\\
\nabla_{h}^{2}\mathbf{A}_{g+1/2}^{n+1/2}+\mathbf{\bar{j}}_{g+1/2}^{n+1/2} & = & \frac{\nabla_{h}\phi_{g}^{n+1}-\nabla_{h}\phi_{g}^{n}}{\Delta t},\label{eq:ampere-md-em}\\
\mathbf{E}_{g+1/2}^{n+1/2} & = & -\nabla_{h}\phi_{g}^{n+1/2}-\frac{\mathbf{A}_{g+1/2}^{n+1}-\mathbf{A}_{g+1/2}^{n}}{\Delta t},\label{eq:e-field-def-em}\\
\mathbf{\bar{j}}_{g+1/2}^{n+1/2} & = & \frac{1}{\Delta_{h}\Delta t}\sum_{p}q_{p}\sum_{\nu\in(n,n+1)}\Delta\tau^{\nu}\mathbf{v}_{p}^{\nu+1/2}\cdot\bar{\bar{\mathbf{S}}}(\mathbf{x}_{g+1/2}-\mathbf{x}_{p}^{\nu+1/2}),\label{eq:j_gather-md-1}\\
\nabla_{h}\cdot\mathbf{A}_{g+1/2}^{n+1} & = & 0,
\end{eqnarray}
where, as before, $\phi_{g}$ is defined at cell centers, and $\mathbf{E}_{g+1/2}$,
$\bar{\mathbf{j}}_{g+1/2}$, and $\mathbf{A}_{g+1/2}$ are defined
at faces.

\subsection{Kinetic energy local conservation law}

The orbit-averaged kinetic-energy update from $n$ to $n+1$ remains
unchanged from the previous analysis: 
\begin{eqnarray*}
\frac{e_{k,g}^{n+1}-e_{k,g}^{n}}{\Delta t}+\nabla_{h}\cdot\boldsymbol{\Gamma}_{g+1/2}^{n+1/2} & = & S_{g}^{n+1/2},
\end{eqnarray*}
with $e_{k}$ defined as before, and fluxes and source defined in
Eqs. \ref{eq:md-source}-\ref{eq:md-flux-z}.

\subsection{Electrostatic energy local conservation law}

The electrostatic field energy conservation law can be derived as
before by considering the local definition: 
\[
\mathcal{E}_{g}=\frac{(E_{x,i+1/2,j,k}^{ES})^{2}+(E_{x,i-1/2,j,k}^{ES})^{2}+(E_{y,i,j+1/2,k}^{ES})^{2}+(E_{y,i,j-1/2,k}^{ES})^{2}+(E_{z,i,j,k+1/2}^{ES})^{2}+(E_{z,i,j,k-1/2}^{ES})^{2}}{4},
\]
where $\mathbf{E}_{g+1/2}^{ES}=-\nabla_{h}\phi_{g}$. The temporal
update equation reads (following the exact same procedure as in the
last section): 
\begin{eqnarray}
\frac{\mathcal{E}_{g}^{n+1}-\mathcal{E}_{g}^{n}}{\Delta t} & = & \nabla_{h}\cdot\boldsymbol{\Gamma}_{E,g+1/2}^{n+1/2}+\frac{1}{2}\left[\bar{\mathbf{j}}_{g+1/2}^{n+1/2}\cdot\nabla_{h}\phi_{g}^{n+1/2}+\bar{\mathbf{j}}_{g-1/2}^{n+1/2}\cdot\nabla_{h}\phi_{g-1}^{n+1/2}\right],\label{eq:es_local_energy}
\end{eqnarray}
where, from Eqs. \ref{eq:E-energy-flux} and \ref{eq:ampere-md-em},
we have: 
\[
\boldsymbol{\Gamma}_{E,g+1/2}^{n+1/2}=-\phi_{g+1/2}^{n+1/2}\left(\frac{\mathbf{E}_{g+1/2}^{ES,n+1}-\mathbf{E}_{g+1/2}^{ES,n}}{\Delta t}+\bar{\mathbf{j}}_{g+1/2}^{n+1/2}\right)=\phi_{g+1/2}^{n+1/2}\nabla_{h}^{2}\mathbf{A}_{g+1/2}^{n+1/2},
\]
which is the discrete form of the first term in the continuum EM flux
(Eq. \ref{eq:em-flux}).

For the second term in Eq. \ref{eq:es_local_energy}, we have, using
Eq. \ref{eq:e-field-def-em}: 
\begin{eqnarray*}
 &  & \frac{1}{2}\left[\bar{\mathbf{j}}_{g+1/2}^{n+1/2}\cdot\nabla_{h}\phi_{g}^{n+1/2}+\bar{\mathbf{j}}_{g-1/2}^{n+1/2}\cdot\nabla_{h}\phi_{g-1}^{n+1/2}\right]=-\frac{\mathbf{E}_{g+1/2}^{n+1/2}\cdot\bar{\mathbf{j}}_{g+1/2}^{n+1/2}+\mathbf{E}_{g-1/2}^{n+1/2}\cdot\bar{\mathbf{j}}_{g-1/2}^{n+1/2}}{2}\\
 & - & \frac{1}{2}\left[\frac{\mathbf{A}_{g+1/2}^{n+1}-\mathbf{A}_{g+1/2}^{n}}{\Delta t}\cdot\bar{\mathbf{j}}_{g+1/2}^{n+1/2}+\frac{\mathbf{A}_{g-1/2}^{n+1}-\mathbf{A}_{g-1/2}^{n}}{\Delta t}\cdot\bar{\mathbf{j}}_{g-1/2}^{n+1/2}\right].
\end{eqnarray*}
Putting things together, we find from Eq. \ref{eq:es_local_energy}:
\begin{eqnarray}
\frac{\mathcal{E}_{g}^{n+1}-\mathcal{E}_{g}^{n}}{\Delta t} & = & \nabla_{h}\cdot\boldsymbol{\Gamma}_{E,g+1/2}^{n+1/2}-\frac{\mathbf{E}_{g+1/2}^{n+1/2}\cdot\bar{\mathbf{j}}_{g+1/2}^{n+1/2}+\mathbf{E}_{g-1/2}^{n+1/2}\cdot\bar{\mathbf{j}}_{g-1/2}^{n+1/2}}{2}\nonumber \\
 & - & \frac{1}{2}\left[\frac{\mathbf{A}_{g+1/2}^{n+1}-\mathbf{A}_{g+1/2}^{n}}{\Delta t}\cdot\bar{\mathbf{j}}_{g+1/2}^{n+1/2}+\frac{\mathbf{A}_{g-1/2}^{n+1}-\mathbf{A}_{g-1/2}^{n}}{\Delta t}\cdot\bar{\mathbf{j}}_{g-1/2}^{n+1/2}\right].\label{eq:es_local_energy_final}
\end{eqnarray}
The last term in this equation encodes the magnetic energy term. We
discuss this next.

\subsection{Magnetic energy local conservation law}

Focusing on the last term in Eq. \ref{eq:es_local_energy_final},
and using Eq. \ref{eq:ampere-md-em} again: 
\begin{eqnarray}
 & - & \frac{1}{2}\left[\frac{\mathbf{A}_{g+1/2}^{n+1}-\mathbf{A}_{g+1/2}^{n}}{\Delta t}\cdot\bar{\mathbf{j}}_{g+1/2}^{n+1/2}+\frac{\mathbf{A}_{g-1/2}^{n+1}-\mathbf{A}_{g-1/2}^{n}}{\Delta t}\cdot\bar{\mathbf{j}}_{g-1/2}^{n+1/2}\right]\nonumber \\
 & = & \frac{1}{2}\left[\frac{\mathbf{A}_{g+1/2}^{n+1}-\mathbf{A}_{g+1/2}^{n}}{\Delta t}\cdot\nabla_{h}^{2}\mathbf{A}_{g+1/2}^{n+1/2}+\frac{\mathbf{A}_{g-1/2}^{n+1}-\mathbf{A}_{g-1/2}^{n}}{\Delta t}\cdot\nabla_{h}^{2}\mathbf{A}_{g-1/2}^{n+1/2}\right]\nonumber \\
 & - & \frac{1}{2}\left[\frac{\mathbf{A}_{g+1/2}^{n+1}-\mathbf{A}_{g+1/2}^{n}}{\Delta t}\cdot\nabla_{h}\frac{\phi_{g}^{n+1}-\phi_{g}^{n}}{\Delta t}+\frac{\mathbf{A}_{g-1/2}^{n+1}-\mathbf{A}_{g-1/2}^{n}}{\Delta t}\cdot\nabla_{h}\frac{\phi_{g-1}^{n+1}-\phi_{g-1}^{n}}{\Delta t}\right].\label{eq:toward_M_energy}
\end{eqnarray}
We begin with the first term in the left-hand-side of this equation.
Discretizing the vector Laplacian in 2D Cartesian geometry for simplicity,
calling $\mathbf{C}_{g+1/2}=(\mathbf{A}_{g+1/2}^{n+1}-\mathbf{A}_{g+1/2}^{n})/{\Delta t}$,
and using Eq. \ref{eq:exact-chain-rule-div}, we can write (with a
sum implied over components $x,y,z$): 
\begin{eqnarray}
 &  & \frac{1}{2}\left[C_{i+1/2,j}^{(x,y,z)}\frac{\partial_{x}A_{i+1,j}^{(x,y,z)}-\partial_{x}A_{i,j}^{(x,y,z)}}{\Delta x}+C_{i-1/2,j}^{(x,y,z)}\frac{\partial_{x}A_{i,j}^{(x,y,z)}-\partial_{x}A_{i-1,j}^{(x,y,z)}}{\Delta x}\right]\nonumber \\
 & + & \frac{1}{2}\left[C_{i,j+1/2}^{(x,y,z)}\frac{\partial_{y}A_{i,j+1}^{(x,y,z)}-\partial_{y}A_{i,j}^{(x,y,z)}}{\Delta y}+C_{i,j-1/2}^{(x,y,z)}\frac{\partial_{y}A_{i,j}^{(x,y,z)}-\partial_{y}A_{i,j-1}^{(x,y,z)}}{\Delta y}\right]\nonumber \\
 & = & \frac{1}{\Delta x}\left[\frac{\partial_{x}A_{i+1,j}^{(x,y,z)}+\partial_{x}A_{i,j}^{(x,y,z)}}{2}C_{i+1/2,j}^{(x,y,z)}-\frac{\partial_{x}A_{i,j}^{(x,y,z)}+\partial_{x}A_{i-1,j}^{(x,y,z)}}{2}C_{i-1/2,j}^{(x,y,z)}\right]-\partial_{x}A_{i,j}^{(x,y,z)}\frac{C_{i+1/2,j}^{(x,y,z)}-C_{i-1/2,j}^{(x,y,z)}}{\Delta x}\nonumber \\
 & + & \frac{1}{\Delta y}\left[\frac{\partial_{y}A_{i,j+1}^{(x,y,z)}+\partial_{y}A_{i,j}^{(x,y,z)}}{2}C_{i,j+1/2}^{(x,y,z)}-\frac{\partial_{y}A_{i,j}^{(x,y,z)}+\partial_{y}A_{i,j-1}^{(x,y,z)}}{2}C_{i,j-1/2}^{(x,y,z)}\right]-\partial_{y}A_{i,j}^{(x,y,z)}\frac{C_{i,j+1/2}^{(x,y,z)}-C_{i,j-1/2}^{(x,y,z)}}{\Delta y}.\label{eq:tensor-chain-rule}
\end{eqnarray}
Using this result, and recalling that (Eq. \ref{eq:veclap_def}):
\[
\nabla^{2}\mathbf{A}=\nabla\cdot(\nabla\mathbf{A}-\nabla\mathbf{A}^{\dagger}),
\]
the first term in the right hand side of Eq. \ref{eq:toward_M_energy}
can be rewritten as: 
\begin{eqnarray*}
 &  & \frac{1}{2}\left[\frac{\mathbf{A}_{g+1/2}^{n+1}-\mathbf{A}_{g+1/2}^{n}}{\Delta t}\cdot\nabla_{h}^{2}\mathbf{A}_{g+1/2}^{n+1/2}+\frac{\mathbf{A}_{g-1/2}^{n+1}-\mathbf{A}_{g-1/2}^{n}}{\Delta t}\cdot\nabla_{h}^{2}\mathbf{A}_{g-1/2}^{n+1/2}\right]\\
 & = & \nabla_{h}\cdot\left((\nabla_{h}\mathbf{A}_{g}^{n+1/2}-\nabla_{h}\mathbf{A}_{g}^{\dagger,n+1/2})\cdot\frac{\mathbf{A}_{g}^{n+1}-\mathbf{A}_{g}^{n}}{\Delta t}\right)-(\nabla_{h}\mathbf{A}_{g}^{n+1/2}-\nabla_{h}\mathbf{A}_{g}^{\dagger,n+1/2}):\nabla_{h}\frac{\mathbf{A}_{g}^{n+1}-\mathbf{A}_{g}^{n}}{\Delta t}.
\end{eqnarray*}
The last term in this expression can in turn be rewritten as (proof
by expansion): 
\[
(\nabla_{h}\mathbf{A}_{g}^{n+1/2}-\nabla_{h}\mathbf{A}_{g}^{\dagger,n+1/2}):\nabla_{h}\frac{\mathbf{A}_{g}^{n+1}-\mathbf{A}_{g}^{n}}{\Delta t}=\frac{\nabla_{h}\mathbf{A}_{g}^{n+1}:(\nabla_{h}\mathbf{A}_{g}^{n+1}-\nabla_{h}\mathbf{A}_{g}^{\dagger,n+1})-\nabla_{h}\mathbf{A}_{g}^{n}:(\nabla_{h}\mathbf{A}_{g}^{n}-\nabla_{h}\mathbf{A}_{g}^{\dagger,n})}{2\Delta t},
\]
and therefore: 
\[
\frac{1}{2}\left[\frac{\mathbf{A}_{g+1/2}^{n+1}-\mathbf{A}_{g+1/2}^{n}}{\Delta t}\cdot\nabla_{h}^{2}\mathbf{A}_{g+1/2}^{n+1/2}+\frac{\mathbf{A}_{g-1/2}^{n+1}-\mathbf{A}_{g-1/2}^{n}}{\Delta t}\cdot\nabla_{h}^{2}\mathbf{A}_{g-1/2}^{n+1/2}\right]=\left.\nabla_{h}\cdot\boldsymbol{\Gamma}_{M}^{n+1/2}\right|_{g}-\frac{\mathcal{M}_{g}^{n+1}-\mathcal{M}_{g}^{n}}{\Delta t},
\]
with: 
\begin{eqnarray*}
\boldsymbol{\Gamma}_{M}^{n+1/2} & = & (\nabla_{h}\mathbf{A}_{g+1/2}^{n+1/2}-\nabla_{h}\mathbf{A}_{g+1/2}^{\dagger,n+1/2})\cdot\frac{\mathbf{A}_{g+1/2}^{n+1}-\mathbf{A}_{g+1/2}^{n}}{\Delta t},\\
\mathcal{M}_{g}^{n+1} & = & \frac{1}{2}\nabla_{h}\mathbf{A}_{g}^{n+1}:(\nabla_{h}\mathbf{A}_{g}^{n+1}-\nabla_{h}\mathbf{A}_{g}^{\dagger,n+1}).
\end{eqnarray*}
The flux $\boldsymbol{\Gamma}_{M}$ corresponds to the second term
in the continuum energy flux definition, Eq. \ref{eq:em-flux}.

Regarding the second term in the right hand side of Eq. \ref{eq:toward_M_energy},
we can exactly write (by Eq. \ref{eq:exact-chain-rule-div}): 
\begin{eqnarray*}
 &  & -\frac{1}{2}\left[\frac{\mathbf{A}_{g+1/2}^{n+1}-\mathbf{A}_{g+1/2}^{n}}{\Delta t}\cdot\nabla_{h}\frac{\phi_{g+1/2}^{n+1}-\phi_{g+1/2}^{n}}{\Delta t}+\frac{\mathbf{A}_{g-1/2}^{n+1}-\mathbf{A}_{g-1/2}^{n}}{\Delta t}\cdot\nabla_{h}\frac{\phi_{g-1/2}^{n+1}-\phi_{g-1/2}^{n}}{\Delta t}\right]\\
 & = & -\nabla_{h}\cdot\left(\frac{\phi_{g+1/2}^{n+1}-\phi_{g+1/2}^{n}}{\Delta t}\frac{\mathbf{A}_{g+1/2}^{n+1}-\mathbf{A}_{g+1/2}^{n}}{\Delta t}\right)+\frac{\phi_{g}^{n+1}-\phi_{g}^{n}}{\Delta t}\cancelto{0}{\nabla_{h}\cdot\frac{\mathbf{A}_{g+1/2}^{n+1}-\mathbf{A}_{g+1/2}^{n}}{\Delta t}}=\nabla_{h}\cdot\boldsymbol{\Gamma}_{M2,g+1/2}^{n+1/2},
\end{eqnarray*}
where we have used the Coulomb gauge, and defined: 
\[
\boldsymbol{\Gamma}_{M2,g+1/2}^{n+1/2}=-\frac{\phi_{g+1/2}^{n+1}-\phi_{g+1/2}^{n}}{\Delta t}\frac{\mathbf{A}_{g+1/2}^{n+1}-\mathbf{A}_{g+1/2}^{n}}{\Delta t}.
\]
As before, $\phi_{g+1/2}$ is the average of cell-centered values
to the corresponding face. This flux is the discrete form of the third
term in the continuum energy flux equation, Eq. \ref{eq:em-flux}.

\subsection{EM total-energy local conservation law}

Putting things together, we find: 
\begin{eqnarray*}
\frac{\mathcal{M}_{g}^{n+1}-\mathcal{M}_{g}^{n}}{\Delta t}+\frac{\mathcal{E}_{g}^{n+1}-\mathcal{E}_{g}^{n}}{\Delta t} & = & \nabla_{h}\cdot\left(\boldsymbol{\Gamma}_{E,g+1/2}^{n+1/2}+\boldsymbol{\Gamma}_{M,g+1/2}^{n+1/2}+\boldsymbol{\Gamma}_{M2,g+1/2}^{n+1/2}\right)\\
 & - & \frac{\mathbf{E}_{g+1/2}^{n+1/2}\cdot\bar{\mathbf{j}}_{g+1/2}^{n+1/2}+\mathbf{E}_{g-1/2}^{n+1/2}\cdot\bar{\mathbf{j}}_{g-1/2}^{n+1/2}}{2}.
\end{eqnarray*}
This is the EM \emph{local} energy conservation law, which when combined
with the kinetic-energy local conservation law, yields the final total
energy local conservation law sought: 
\begin{eqnarray*}
\frac{e_{k,g}^{n+1}-e_{k,g}^{n}}{\Delta t}+\frac{\mathcal{M}_{g}^{n+1}-\mathcal{M}_{g}^{n}}{\Delta t}+\frac{\mathcal{E}_{g}^{n+1}-\mathcal{E}_{g}^{n}}{\Delta t} & = & \nabla_{h}\cdot\left(-\boldsymbol{\Gamma}_{g+1/2}^{n+1/2}+\boldsymbol{\Gamma}_{E,g+1/2}^{n+1/2}+\boldsymbol{\Gamma}_{M,g+1/2}^{n+1/2}+\boldsymbol{\Gamma}_{M2,g+1/2}^{n+1/2}\right)\\
 & + & \frac{1}{\Delta_{h}\Delta t}\sum_{p}q_{p}\sum_{\nu\in(n,n+1)}\Delta\tau^{\nu}\mathbf{E}_{p}^{\nu+1/2}\cdot\mathbf{v}_{p}^{\nu+1/2}\tilde{\tilde{S}}_{l}(\mathbf{x}_{g}-\mathbf{x}_{p}^{\nu+1/2})\\
 & - & \frac{1}{2\Delta_{h}\Delta t}\sum_{p}q_{p}\sum_{\nu\in(n,n+1)}\Delta\tau^{\nu}\mathbf{v}_{p}^{\nu+1/2}\cdot\bar{\bar{\mathbf{S}}}(\mathbf{x}_{g+1/2}-\mathbf{x}_{p}^{\nu+1/2})\cdot\mathbf{E}_{g+1/2}^{n+1/2}\\
 & - & \frac{1}{2\Delta_{h}\Delta t}\sum_{p}q_{p}\sum_{\nu\in(n,n+1)}\Delta\tau^{\nu}\mathbf{v}_{p}^{\nu+1/2}\cdot\bar{\bar{\mathbf{S}}}(\mathbf{x}_{g-1/2}-\mathbf{x}_{p}^{\nu+1/2})\cdot\mathbf{E}_{g-1/2}^{n+1/2}.
\end{eqnarray*}
The last term vanishes with mesh/particle refinement, and can be written
in a fully conservative form, as before. We therefore finally arrive
at the full total-energy electromagnetic \emph{local} conservation
law: 
\[
\frac{e_{k,g}^{n+1}-e_{k,g}^{n}}{\Delta t}+\frac{\mathcal{M}_{g}^{n+1}-\mathcal{M}_{g}^{n}}{\Delta t}+\frac{\mathcal{E}_{g}^{n+1}-\mathcal{E}_{g}^{n}}{\Delta t}+\nabla_{h}\cdot\left(\boldsymbol{\Gamma}_{g+1/2}^{n+1/2}-\boldsymbol{\Gamma}_{E,g+1/2}^{n+1/2}-\boldsymbol{\Gamma}_{M,g+1/2}^{n+1/2}-\boldsymbol{\Gamma}_{M2,g+1/2}^{n+1/2}\right)=\nabla_{h}\cdot\boldsymbol{\Gamma}_{K,g+1/2}^{n+1/2}.
\]

\section{Numerical results}

\label{sec:numerics} In the previous local-energy conservation analysis
of the various implicit PIC models, a numerical flux $\boldsymbol{\Gamma}_{K}$
appears that remains unquantified. The presence of this unquantified
numerical term is a weakness of the conclusion of the analysis, because
it is not clear \emph{a priori} the relative importance of the contribution
of this term to the local energy balance. 

The objective of this section is two-fold: 1) to quantify the importance
of the numerical flux to the local energy balance in a representative
multiscale problem with reasonable resolution; 2) to verify the key
ingredients of the local-conservation proofs outlined in this study.
For this purpose, we have implemented diagnostics in the fully implicit
electromagnetic code DPIC~\citep{chen2015multi} to quantify all
the contributions to the local energy balance equation, including
the numerical term. We consider the electrostatic formulation of Sec.
\ref{sec:multi-d} with orbit averaging and an imposed magnetic field,
but specialized to one dimension. Multiple dimensions and electromagnetic
effects do not change the form of the numerical energy source, and
involve similar manipulations as in the electrostatic case to prove
the local energy conservation theorem, and therefore are not essential
to the objectives of this section.

\begin{figure}
\centering{}\includegraphics[width=0.8\columnwidth]{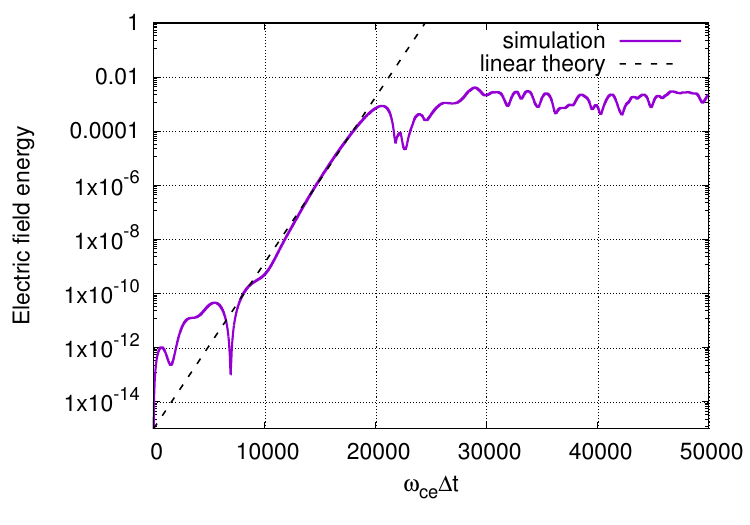}\caption{Modified two-stream instability problem: Electrostatic energy growth
and saturation with timestep $\Delta t=0.2\omega_{ce}^{-1}$.}
\label{fig:modified-2stream} 
\end{figure}

\begin{figure}
\centering{}\includegraphics[width=0.8\columnwidth]{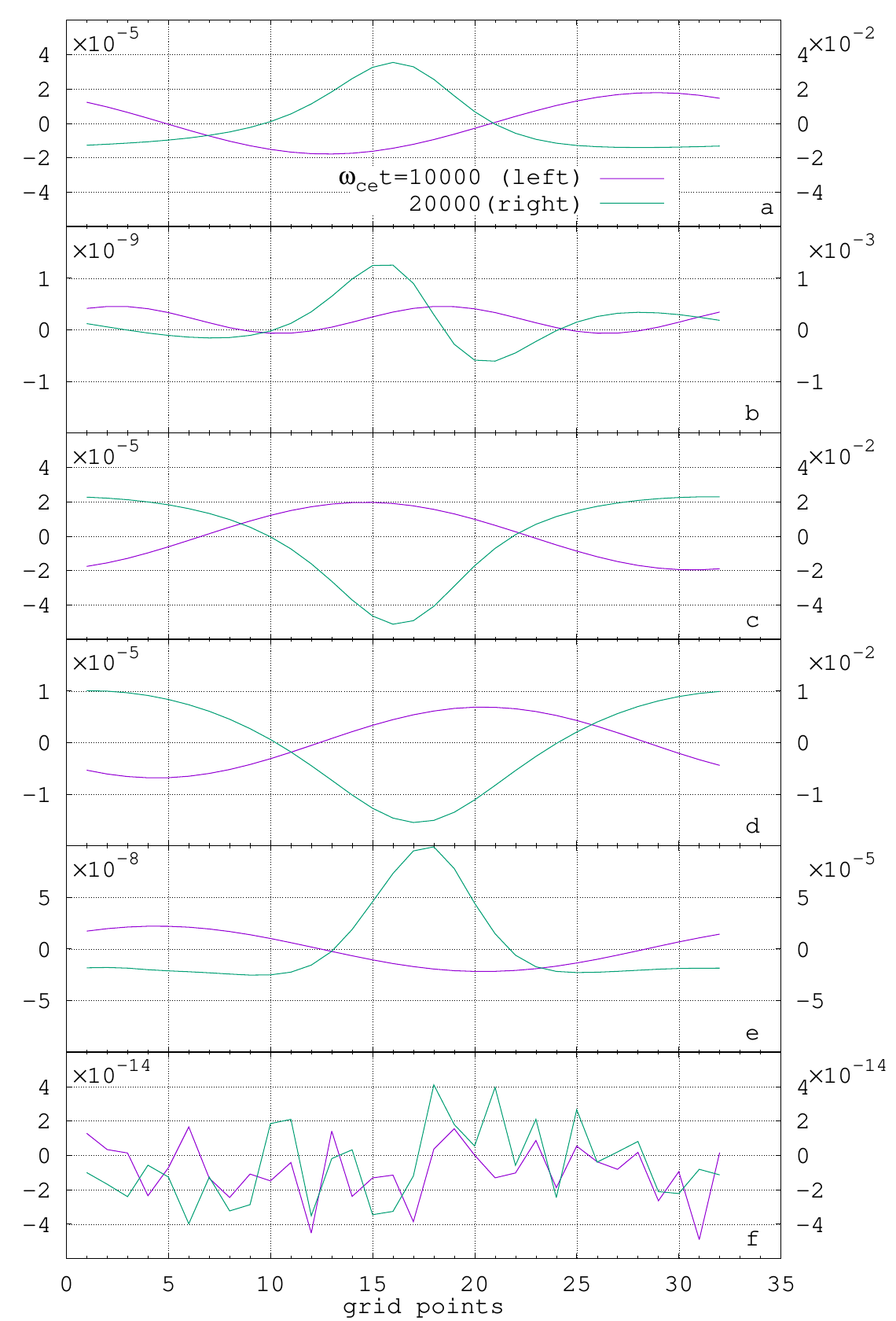}\caption{Modified two-stream instability problem: Simulation results of various
terms in the energy balance equation at two different times: $\omega_{ce}t=10000,20000$.
Panels a-e depict $(e_{k,g}^{n+1}-e_{k,g}^{n})/\Delta t$, $(\mathcal{E}_{g}^{n+1}-\mathcal{E}_{g}^{n})/\Delta t$,
$\nabla_{h}\cdot\boldsymbol{\Gamma}_{g+1/2}$, $\nabla_{h}\cdot\boldsymbol{\Gamma}_{E,g+1/2}$
and $\nabla_{h}\cdot\boldsymbol{\Gamma}_{K,g+1/2}$ in Eq.~\ref{eq:energy-balance-es1d}
respectively. Panel f denotes the residue of the balance equation,
which vanishes to numerical round-off.}
\label{fig:modified-2stream-terms} 
\end{figure}

We consider a modified two-stream instability (MTSI) problem for our
numerical tests, which is a multiscale variant of the classic two-stream
instability. MTSI occurs in a plasma when there is a relative drift
between electrons and ions in the presence of a magnetic field \citep{mcbride1972theory}.
The simulation is performed in 1D-3V with $L_{x}=1.8229$ $\lambda_{Di}$,
$\omega_{ce}/\omega_{pe}=10$ and $m_{i}/m_{e}=5000$. In ion units,
$m_{i}=1.0$ and $B=1/\sqrt{50}$. The magnetic field is mostly pointed
in the $y$ direction but slightly tilted such that $B_{x}/B=\sin(\theta)=\sqrt{m_{e}/m_{i}}$.
Electrons are initially stationary, and ions have a relative velocity
of $U=0.5$. Both species are cold initially. For these parameters,
the growth rate is $\gamma=0.4992\omega_{pi}.$ The base resolution
is $N_{g}=32$ cells, $N_{ppc}=100$ particles per cell, and a timestep
$\Delta t=0.2\omega_{ce}^{-1}$, small enough to resolve the electron
gyromotion. We run the simulation up to $\omega_{ce}t=50000$, long
enough for linear growth and nonlinear saturation of the MTSI. Figure~\ref{fig:modified-2stream}
shows excellent agreement in the linear growth rate between simulation
and linear theory. Figure ~\ref{fig:modified-2stream-terms} depicts
the magnitude of various terms at various times during the simulation
in the electrostatic local energy balance equation: 
\begin{equation}
\frac{\left(e_{k,g}^{n+1}-e_{k,g}^{n}\right)}{\Delta t}+\frac{\left(\mathcal{E}_{g}^{n+1}-\mathcal{E}_{g}^{n}\right)}{\Delta t}+\nabla_{h}\cdot\left(\boldsymbol{\Gamma}_{g+1/2}^{n+1/2}-\boldsymbol{\Gamma}_{E,g+1/2}^{n+1/2}\right)=\nabla_{h}\cdot\boldsymbol{\Gamma}_{K,g+1/2}^{n+1/2}.\label{eq:energy-balance-es1d}
\end{equation}
The figure demonstrates that this equation indeed vanishes to numerical
round off, supporting the analysis. The figure also demonstrates that
the numerical flux term is subdominant by several orders of magnitude
vs. other contributions, particularly in the nonlinear stage.

We test numerical convergence of the error term, $\nabla_{h}\cdot\boldsymbol{\Gamma}_{K,g+1/2}$,
by performing a grid, timestep and particle-number convergence study.
Assuming a second-order method in time and space, and the usual $1/\sqrt{N_{ppc}}$
Monte Carlo scaling of the error due to particle noise, increasing
the spatial mesh by a factor of 2 requires decreasing the timestep
by the same factor and increasing the number of particles per cell
by a factor of 16 to ensure that all error contributions decrease
commensurately by a factor of 4 \citep{tranquilli2022deterministic}.
This results in simulation parameters $N_{g}=64$, $N_{ppc}=1600$,
and $\omega_{ce}\Delta t=0.1$. The expected error decrease is confirmed
in Table \ref{table:1} at two different time snapshots in the simulation,
$\omega_{ce}t=10000,20000$, confirming the proper implementation
of the scheme and the asymptotic cancellation of the numerical error
term. 
\begin{table}[h!]
\centering{}\caption{Convergence test of the error term in the local energy balance equation
at two time snapshots for two simulations with coarse and fine resolutions
as described in the main text. The error is computed as err=$||\nabla_{h}\cdot\boldsymbol{\Gamma}_{K,g+1/2}||_{1}/N_{c}$,
where $N_{c}$ is the number of grid points.}
\label{table:1} %
\begin{tabular}{c|cc}
$\omega_{ce}t$  & 10000  & 20000\tabularnewline
\hline 
err$_{c}$  & $1.4\times10^{-8}$  & $2.92\times10^{-5}$\tabularnewline
err$_{f}$  & $3.54\times10^{-9}$  & $7.57\times10^{-8}$\tabularnewline
\hline 
ratio  & 3.95  & 3.86\tabularnewline
\end{tabular}
\end{table}

\section{Conclusions}

\label{sec:conclusions}

In this study, we have demonstrated rigorously that the class of fully
implicit, local-charge- and global-energy-conserving PIC algorithms
proposed over the last decade \citep{chen2011energy,chen20141ddarwin,chacon2013curvilinear,chen2015multi,chacon2016curvilinear,chen2023implicit}
feature a \emph{local} total-energy conservation theorem as well.
This result evidences the strong theoretical foundations of these
fully implicit PIC algorithms. The proof leverages recent analysis
for geometric PIC algorithms \citep{xiao2017local}; however, unlike
those studies (which were performed with continuous time and were
not numerically verified), our result is time-discrete (including
orbit averaging) and has been verified numerically, albeit in one
dimension for simplicity. Our numerical results verify that a local
energy conservation principle indeed exists, and demonstrate that
the numerical energy flux term present in the local energy balance
equation is subdominant vs. physical fluxes for representative problems.

\section*{Acknowledgments}

The authors would like to acknowledge useful conversations with J.
Angus and H. Qin, and in particular for their prodding us to look
in this direction. This work was sponsored by the DOE Office of Applied
Scientific Computing Research Mathematical Multifaceted Integrated
Capability Centers (MMICCs) and Applied Mathematics Research programs.
This work was performed at Los Alamos National Laboratory, operated
by Triad National Security, LLC, for the National Nuclear Security
Administration of the U.S. Department of Energy (Contract No. 89233218CNA000001).

\appendix

\section{Derivation of local energy conservation law in 1D electrostatic model}

\label{app:1d-es-derivation}

Starting from Eq. \ref{eq:dEk_dt} and introducing Eqs. , \ref{eq:taylor-n}\ref{eq:taylor-np},
we find: 
\begin{eqnarray*}
\frac{e_{k,i}^{n+1}-e_{k,i}^{n}}{\Delta t} & = & \frac{1}{\Delta x\Delta t}\sum_{p}m_{p}\left[\frac{(v_{p}^{n+1})^{2}}{2}-\frac{(v_{p}^{n})^{2}}{2}\right]S_{l}(x_{i}-x_{p}^{n+1/2})\\
 & - & \frac{1}{\Delta x\Delta t}\sum_{p}m_{p}\left.\frac{\partial S_{l}}{\partial x_{p}}\right|_{i}^{n+1/2}\left[\frac{(v_{p}^{n+1})^{2}}{2}(x_{p}^{n+1}-x_{p}^{n+1/2})-\frac{(v_{p}^{n})^{2}}{2}(x_{p}^{n}-x_{p}^{n+1/2})\right]\\
 & + & \frac{1}{\Delta x\Delta t}\sum_{p}m_{p}\left.\frac{\partial^{2}S_{l}}{\partial x_{p}^{2}}\right|_{i}^{n+1/2}\left[\frac{(v_{p}^{n+1})^{2}}{2}\frac{(x_{p}^{n+1}-x_{p}^{n+1/2})^{2}}{2}-\frac{(v_{p}^{n})^{2}}{2}\frac{(x_{p}^{n}-x_{p}^{n+1/2})^{2}}{2}\right].
\end{eqnarray*}
Working on the terms of each sum independently, the first term reads
(using Eq. \ref{eq:eom}): 
\[
\frac{m_{p}}{\Delta t}\left[\frac{(v_{p}^{n+1})^{2}}{2}-\frac{(v_{p}^{n})^{2}}{2}\right]=m_{p}\frac{(v_{p}^{n+1}-v_{p}^{n})}{\Delta t}v_{p}^{n+1/2}=q_{p}E_{p}^{n+1/2}v_{p}^{n+1/2}.
\]
The second term reads, using that $(x_{p}^{n+1}-x_{p}^{n+1/2})=\Delta tv_{p}^{n+1/2}/2$
by definition of $x_{p}^{n+1/2}$ (and similarly with other terms):
\[
\frac{m_{p}}{\Delta x\Delta t}\left.\frac{\partial S_{l}}{\partial x_{p}}\right|_{i}^{n+1/2}\left[\frac{(v_{p}^{n+1})^{2}}{2}(x_{p}^{n+1}-x_{p}^{n+1/2})-\frac{(v_{p}^{n})^{2}}{2}(x_{p}^{n}-x_{p}^{n+1/2})\right]=\frac{1}{\Delta x}\left.\frac{\partial S_{l}}{\partial x_{p}}\right|_{i}^{n+1/2}v_{p}^{n+1/2}e_{k,p}^{n+1/2},
\]
where: 
\[
e_{k,p}^{n+1/2}=\frac{m_{p}}{2}\left[\frac{(v_{p}^{n+1})^{2}}{2}+\frac{(v_{p}^{n})^{2}}{2}\right]
\]
is a time-centered particle kinetic energy. Noting that: 
\[
\left.\frac{\partial S_{l}}{\partial x_{p}}\right|_{i}^{n+1/2}=\frac{1}{\Delta x}\left[S_{l-1}(x_{i+1/2}-x_{p}^{n+1/2})-S_{l-1}(x_{i-1/2}-x_{p}^{n+1/2})\right],
\]
we find: 
\begin{eqnarray*}
\frac{m_{p}}{\Delta x\Delta t}\left.\frac{\partial S_{l}}{\partial x_{p}}\right|_{i}^{n+1/2}\left[\frac{(v_{p}^{n+1})^{2}}{2}(x_{p}^{n+1}-x_{p}^{n+1/2})-\frac{(v_{p}^{n})^{2}}{2}(x_{p}^{n}-x_{p}^{n+1/2})\right]\\
=\frac{S_{l-1}(x_{i+1/2}-x_{p}^{n+1/2})-S_{l-1}(x_{i-1/2}-x_{p}^{n+1/2})}{\Delta x}v_{p}^{n+1/2}\frac{e_{k,p}^{n+1/2}}{\Delta x} & ,
\end{eqnarray*}
from which a flux form can be derived.

Finally, the third term gives: 
\begin{eqnarray*}
\frac{(v_{p}^{n+1})^{2}}{2}\frac{(x_{p}^{n+1}-x_{p}^{n+1/2})^{2}}{2}-\frac{(v_{p}^{n})^{2}}{2}\frac{(x_{p}^{n}-x_{p}^{n+1/2})^{2}}{2} & = & \frac{\Delta t^{2}}{8}\left(v_{p}^{n+1/2}\right)^{2}\left(\frac{(v_{p}^{n+1})^{2}}{2}-\frac{(v_{p}^{n})^{2}}{2}\right)\\
=\frac{\Delta t^{2}}{8}\left(v_{p}^{n+1/2}\right)^{3}(v_{p}^{n+1}-v_{p}^{n}) & = & \frac{\Delta t^{3}}{8}\left(v_{p}^{n+1/2}\right)^{3}\frac{q_{p}}{m_{p}}E_{p}^{n+1/2},
\end{eqnarray*}
which is a second-order temporal correction to the $E.j$ source (and
therefore asymptotically consistent). Putting all together, we find:
\begin{eqnarray}
\frac{e_{k,i}^{n+1}-e_{k,i}^{n}}{\Delta t} & = & \frac{1}{\Delta x}\sum_{p}q_{p}E_{p}^{n+1/2}v_{p}^{n+1/2}\left[S_{l}(x_{i}-x_{p}^{n+1/2})+\left.\frac{\partial^{2}S_{l}}{\partial x_{p}^{2}}\right|_{i}^{n+1/2}\frac{\Delta t^{2}}{8}\left(v_{p}^{n+1/2}\right)^{2}\right]\nonumber \\
 & - & \frac{1}{\Delta x}\sum_{p}\frac{S_{l-1}(x_{i+1/2}-x_{p}^{n+1/2})-S_{l-1}(x_{i-1/2}-x_{p}^{n+1/2})}{\Delta x}v_{p}^{n+1/2}e_{k,p}^{n+1/2}.\label{eq:cons-e_k-2}
\end{eqnarray}
We define: 
\[
\tilde{S}_{l}(x_{i}-x_{p}^{n+1/2})=S_{l}(x_{i}-x_{p}^{n+1/2})+\left.\frac{\partial^{2}S_{l}}{\partial x_{p}^{2}}\right|_{i}^{n+1/2}\frac{\Delta t^{2}}{8}\left(v_{p}^{n+1/2}\right)^{2},
\]
which is also a partition of unity, since: 
\begin{equation}
\sum_{i}\tilde{S}_{l}(x_{i}-x_{p}^{n+1/2})=1,\label{eq:partition-of-unity-1}
\end{equation}
because $\sum_{i}\left.\frac{\partial^{2}S_{l}}{\partial x_{p}^{2}}\right|_{i}^{n+1/2}=0$.
The energy fluxes are defined as: 
\[
\Gamma_{i+1/2}^{n+1/2}=\frac{1}{\Delta x}\sum_{p}S_{l-1}(x_{i+1/2}-x_{p}^{n+1/2})v_{p}^{n+1/2}e_{k,p}^{n+1/2},
\]
and give the following kinetic-energy \emph{local} conservation law:
\begin{equation}
\frac{e_{k,i}^{n+1}-e_{k,i}^{n}}{\Delta t}+\frac{\Gamma_{i+1/2}^{n+1/2}-\Gamma_{i-1/2}^{n+1/2}}{\Delta x}=\frac{1}{\Delta x}\sum_{p}q_{p}E_{p}^{n+1/2}v_{p}^{n+1/2}\tilde{S}_{l}(x_{i}-x_{p}^{n+1/2}),\label{eq:cons-e_k-final-1}
\end{equation}
which is Eq. \ref{eq:cons-e_k-final} in the main text.

\section{Derivation of local energy conservation law in multi-D electrostatic
model}

\label{app:=00003D000020multi-D-es-derivation}

Starting from the orbit-averaged kinetic-energy update from $n$ to
$n+1$ (from Eq. \ref{eq:local-kin-en-1}): 
\begin{equation}
\frac{e_{k,g}^{n+1}-e_{k,g}^{n}}{\Delta t}=\frac{1}{\Delta_{h}\Delta t}\sum_{p}\sum_{\nu\in(n,n+1)}m_{p}\left[\frac{(v_{p}^{\nu+1})^{2}}{2}S_{l}(\mathbf{x}_{g}-\mathbf{x}_{p}^{\nu+1})-\frac{(v_{p}^{\nu})^{2}}{2}S_{l}(\mathbf{x}_{g}-\mathbf{x}_{p}^{\nu})\right].\label{eq:dEk_dt-multi-D-1}
\end{equation}
Taylor-expanding the tensor-product shape function, we have: 
\begin{eqnarray*}
S_{l}(\mathbf{x}_{g}-\mathbf{x}_{p}^{\nu+1}) & = & S_{l}(\mathbf{x}_{g}-\mathbf{x}_{p}^{\nu+1/2})-(\mathbf{x}_{p}^{\nu+1}-\mathbf{x}_{p}^{\nu+1/2})\cdot\left.\nabla_{x_{p}}S_{l}\right|_{g}^{\nu+1/2}\\
 & + & \left.\nabla_{x_{p}}\nabla_{x_{p}}S_{l}\right|_{g}^{\nu+1/2}:\frac{(\mathbf{x}_{p}^{\nu+1}-\mathbf{x}_{p}^{\nu+1/2})\otimes(\mathbf{x}_{p}^{\nu+1}-\mathbf{x}_{p}^{\nu+1/2})}{2},\\
S_{l}(\mathbf{x}_{g}-\mathbf{x}_{p}^{\nu}) & = & S_{l}(\mathbf{x}_{g}-\mathbf{x}_{p}^{\nu+1/2})-(\mathbf{x}_{p}^{\nu}-\mathbf{x}_{p}^{\nu+1/2})\cdot\left.\nabla_{x_{p}}S_{l}\right|_{g}^{\nu+1/2}\\
 & + & \left.\nabla_{x_{p}}\nabla_{x_{p}}S_{l}\right|_{g}^{\nu+1/2}:\frac{(\mathbf{x}_{p}^{\nu}-\mathbf{x}_{p}^{\nu+1/2})\otimes(\mathbf{x}_{p}^{\nu}-\mathbf{x}_{p}^{\nu+1/2})}{2}.
\end{eqnarray*}
There results: 
\begin{eqnarray}
\frac{e_{k,g}^{n+1}-e_{k,g}^{n}}{\Delta t} & = & \frac{1}{\Delta_{h}\Delta t}\sum_{p}\sum_{\nu\in(n,n+1)}\Delta\tau^{\nu}q_{p}\mathbf{E}_{p}^{\nu+1/2}\cdot\mathbf{v}_{p}^{\nu+1/2}\tilde{\tilde{S}}_{l}(\mathbf{x}_{g}-\mathbf{x}_{p}^{\nu+1/2})\nonumber \\
 & - & \frac{1}{\Delta_{h}\Delta t}\sum_{p}\sum_{\nu\in(n,n+1)}\Delta\tau^{\nu}\mathbf{v}_{p}^{\nu+1/2}\cdot\left.\nabla_{x_{p}}S_{l}\right|_{g}^{\nu+1/2}e_{k,p}^{\nu+1/2}\nonumber \\
 & = & S_{g}^{n+1/2}-\nabla_{h}\cdot\boldsymbol{\Gamma}_{g}^{n+1/2},\label{eq:cons-e_k-multiD-1}
\end{eqnarray}
where: 
\[
\tilde{\tilde{S}}_{l}(\mathbf{x}_{g}-\mathbf{x}_{p}^{\nu+1/2})=S_{l}(\mathbf{x}_{g}-\mathbf{x}_{p}^{\nu+1/2})+\frac{(\Delta\tau^{\nu})^{2}}{8}\mathbf{v}_{p}^{\nu+1/2}\cdot(\left.\nabla_{x_{p}}\nabla_{x_{p}}S_{l}\right|_{g}^{\nu+1/2})\cdot\mathbf{v}_{p}^{\nu+1/2},
\]
and, as before: 
\[
e_{k,p}^{\nu+1/2}=\frac{m_{p}}{2}\left[\frac{(v_{p}^{\nu+1})^{2}}{2}+\frac{(v_{p}^{\nu})^{2}}{2}\right].
\]
The source $S_{g}^{n+1/2}$ and the components of the flux $\boldsymbol{\Gamma}_{g+1/2}^{n+1/2}$
are given as: 
\begin{eqnarray*}
S_{g}^{n+1/2} & = & \frac{1}{\Delta_{h}\Delta t}\sum_{p}\sum_{\nu\in(n,n+1)}\Delta\tau^{\nu}q_{p}\mathbf{E}_{p}^{\nu+1/2}\cdot\mathbf{v}_{p}^{\nu+1/2}\tilde{\tilde{S}}_{l}(\mathbf{x}_{g}-\mathbf{x}_{p}^{\nu+1/2}),\\
\Gamma_{x,i+1/2,j,k}^{\nu+1/2} & = & \frac{1}{\Delta_{h}\Delta t}\sum_{p}\sum_{\nu\in(n,n+1)}\Delta\tau^{\nu}S_{1}(x_{i+1/2}-x_{p}^{\nu+1/2})S_{2}(y_{j}-y_{p}^{\nu+1/2})S_{2}(z_{k}-z_{p}^{\nu+1/2})v_{x,p}^{\nu+1/2}e_{k,p}^{\nu+1/2},\\
\Gamma_{y,i,j+1/2,k}^{\nu+1/2} & = & \frac{1}{\Delta_{h}\Delta t}\sum_{p}\sum_{\nu\in(n,n+1)}\Delta\tau^{\nu}S_{2}(x_{i}-x_{p}^{\nu+1/2})S_{1}(y_{j+1/2}-y_{p}^{\nu+1/2})S_{2}(z_{k}-z_{p}^{\nu+1/2})v_{y,p}^{\nu+1/2}e_{k,p}^{\nu+1/2},\\
\Gamma_{z,i,j,k+1/2}^{\nu+1/2} & = & \frac{1}{\Delta_{h}\Delta t}\sum_{p}\sum_{\nu\in(n,n+1)}\Delta\tau^{\nu}S_{2}(x_{i}-x_{p}^{\nu+1/2})S_{2}(y_{j}-y_{p}^{\nu+1/2})S_{1}(z_{k+1/2}-z_{p}^{\nu+1/2})v_{z,p}^{\nu+1/2}e_{k,p}^{\nu+1/2}.
\end{eqnarray*}

 \bibliographystyle{IEEEtran}
\bibliography{refs}

\end{document}